\documentclass{amsart}

\usepackage{amsfonts,amssymb}
\usepackage{comment}

\def\x{\chi}
\def\<{\langle}
\def\>{\rangle}
\def\pf{\noindent {\it Proof.\ }}
\def\ar{\rightarrow}

\def\aut{\rm Aut}

\def\g{\gamma}
\def\Z{\mathbb Z}
\def\gr{\rm gr}
\newcommand{\la}{\lambda}
\newcommand{\ka}{\kappa}

\title{IA-automorphisms and Lie Algebras related to the McCool group}
\author{C.E. Kofinas, V. Metaftsis, A.I. Papistas}
\address{Department of Mathematics, Aristotle University of Thessaloniki, GR-54124 Thessaloniki, Greece.}\email{apapist@math.auth.gr}
\address{Department of Mathematics, University of the Aegean, Karlovassi GR-83200, Samos, Greece}\email{vmet@aegean.gr}
\address{Department of Mathematics, Aristotle University of Thessaloniki, GR-54124 Thessaloniki, Greece}\email{kkofinas@hotmail.com}

\newtheorem{theorem}{Theorem}

\newtheorem{proposition}{Proposition}
\newtheorem{cor}{Corollary}
\newtheorem{corollary}{Corollary}
\newtheorem{lemma}{Lemma}
\newtheorem{remark}{Remark}

\begin{document}

\begin{abstract}
We investigate a subgroup $I_n$ of the McCool group $M_n$. We show that $I_n$ has solvable conjugacy problem. Next, we investigate its Lie algebra $\gr(I_n)$ and we find a presentation for it. Finally, we show that  $\gr(I_n)$ is naturally embedded into the Andreadakis-Johnson Lie algebra of the IA-automorphisms of the free group $F_n$. 
\end{abstract}

\maketitle

\section{Introduction and notation}

For a positive integer $n$, with $n \geq 2$, let $F_{n}$ be a free group of rank $n$ generated by the set $\{ x_1,\ldots, x_n\}$ and $M_n$ be the McCool group (the basis-conjugating group), that is, the subgroup of Aut$(F_{n})$ generated by the automorphism $\x_{ij}$, $1 \leq i\neq j \leq n$, 
with $$\x_{ij}(x_k)=\left\{\begin{array}{lr} x_k & \mbox{if}\ \ k\neq i\\ x^{-1}_jx_ix_j &    \mbox{if}\ \ k=i. \end{array}\right.$$ For any group $G$ and $a, b \in G$, we write $[a,b] = a^{-1}b^{-1}ab$. As shown by McCool in \cite{mccool}, a complete set of relations for the group $M_n$ is the following
$$[\x_{ij},\x_{kj}]=[\x_{ij},\x_{kl}]=[\x_{ij}\x_{kj},\x_{ik}]=1,$$ where different letters correspond to different numbers. It is shown by Savushkina \cite{sav} (see, also, \cite[Theorem 1]{bardakov} and \cite[Theorem 1.2]{cohen})  that $M_n$ is structured as a semi-direct product $M_n=D_{n-1} \rtimes M_{n-1}$, where $D_{n-1}$ is the group generated by the set $\{\x_{1n},\ldots,\x_{n-1,n},\x_{n1},\ldots,$ $\x_{n,n-1}\}.$ The group structure of $D_n$ at the moment is a mystery,  although we know that the subgroup of $D_{n}$ generated by the set $\{\x_{1n},\ldots,\x_{n-1,n}\}$ is free abelian  and the subgroup of $D_{n}$ generated by the set $\{\x_{n1},\ldots,\x_{n,n-1}\}$ is free.  

In the present work, we investigate a subgroup of $M_n$, that arose in \cite{barnes}, the partial inner automorphisms group $I_n$, consisting of automorphisms acting on initial segments of the basis of $F_n$ by conjugations. In Section 2, we show that the conjugacy problem for $I_n$ is solvable, thus answering a question of Bardakov and Neshchadim in \cite{barnes}. In Section 3, we develop an inductive technique that helps us to decompose and understand the graded Lie algebra $\gr(I_n)$ of $I_n$ which leads us to give a presentation for $\gr(I_n)$. Finally, in Section 4, by using standard arguments, we prove that $\gr(I_n)$ naturally embeds into the Andreadakis-Johnson Lie algebra of the IA-automorphisms of the free group $F_n$.

\section{Some group theoretic results} 
 
Bardakov and Neshchadim in \cite{barnes} have shown that in $M_n$, we may define the subgroup $H_n$ generated by the elements $y_{ni}=\x_{1i}\ldots\x_{ni}$, $i=1,\ldots,n$, with the assumption that $\x_{ii}=1$. It is easy to see that $y_{ni}$ is an inner automorphism of $F_{n}$, conjugating each generator by $x_i$. Obviously, the subgroup  $H_n$ of $M_{n}$ generated by the set $\{y_{ni}, i=1,\ldots,n\}$ is equal to Inn$(F_{n})$, the inner automorphisms of $F_{n}$. Similarly, for $k = 2, \ldots, n-1$, they define $H_k$ to be the subgroup of $M_{n}$ generated by the set $\{y_{k1},\ldots, y_{kk}\}$, the inner automorphisms of $F_k$ considered as a subgroup of $M_n$. Finally, they define $I_n$ to be the subgroup of $M_{n}$ generated by the set $\{y_{\mu j}: 2 \leq \mu \leq n, 1 \leq j \leq \mu\}$. Now, one can easily verify that for $i>j$ the group $H_j$ normalizes $H_i$, namely,  
since $H_i$ is generated by the set $\{y_{i1},\ldots, y_{ii}\}$ and $H_j$ is generated by the set $\{y_{j1},\ldots, y_{jj}\}$,  we have
$$y_{jk}^{-1}y_{il}y_{jk}=\left\{\begin{array}{lr} y_{il} & k=l\\y_{il} & l>j\\y^{-1}_{ik}y_{il}y_{ik} & k\neq l, l\le j.\end{array}\right.$$ Furthermore, it is easily shown that $I_n$ is a poly-free group (see also \cite{barnes}). Indeed, if $G_i=H_n\rtimes H_{n-1}\rtimes \ldots\rtimes H_i$, with $i \in \{2, \ldots, n\}$, then there is the tower of normal subgroups $$\{1\} = G_{n+1} \lhd G_n\lhd G_{n-1}\lhd\ldots\lhd G_2=I_n$$ such that each $G_i/G_{i+1}\cong H_{i}$ is free. Thus, $I_{n} = G_{3} \rtimes H_{2}$ and so, any element of $I_{n}$ is uniquely written as $g_{3}w_{2}$,  where $g_{3} \in G_{3}$ and $w_{2} \in H_{2}$. It is shown in \cite{barnes} that a complete set of defining relations for $I_n$ is the following
\begin{enumerate}
\item $[y_{mi},y_{ri}] = 1, ~~~~~~2 \leq r < m \leq n,$
\item $[y_{mi},y_{rj}] = 1, \ \  2 \leq r < i \leq m \leq n, ~~1 \leq j \leq r,$
\item $[y_{mi},y_{rj}]=[y_{mi},y_{mj}],\ \ 2 \leq r < m \leq n,~~ 1 \leq i \leq r < m \leq n,~~ 1 \leq j \leq r, j \neq i.$
\end{enumerate}
Since the group $H_j$ normalizes $H_i$ for $i>j$, the elements of $I_n$ have a normal form. That is, every element of $I_n$ may be uniquely written in the form $w_nw_{n-1}\ldots w_2$, where each $w_i\in H_i = {\rm Inn}(F_{i}) \cong F_i$. This implies that the word problem in $I_n$ is solvable. Indeed, an element $w_nw_{n-1}\ldots w_2=1$ if and only if each $w_i=1$ and that is decidable in a free group. Moreover, the above allow us to answer positively a question of Bardakov and Neshchadim in \cite{barnes}.

\begin{theorem}
The conjugacy problem in $I_n$ is solvable.
\end{theorem}

\pf Let $x=w_n\ldots w_3w_2$, $y=z_n\ldots z_3z_2$ and $g=g_n\ldots g_3g_2$ $\in I_n$, where $w_{i}, z_{i}, g_{i} \in H_{i}$, $i = 2, \ldots, n$, such that $gxg^{-1}=y$. The equation can be rewritten as 
$$g_ng_{n-1}\ldots g_3g_2  w_nw_{n-1}\ldots w_3w_2  g_2^{-1}g_3^{-1}\ldots g_{n-1}^{-1}g_n^{-1}=y$$
or equivalently,
$$g_n  (g_{n-1}\ldots g_2w_ng_2^{-1}\ldots g_{n-1}^{-1})  g_{n-1}  (g_{n-2}\ldots g_3g_2 w_{n-1} g_2^{-1}\ldots g_{n-2}^{-1}) \ldots $$
$$g_4  (g_3g_2w_4g_2^{-1}g_3^{-1})  g_3 (g_2w_3g_2^{-1})(g_{2}w_{2}g^{-1}_{2}) g_{3}^{-1}\ldots g_n^{-1} = y.$$
Write $$a_n=g_{n-1}\ldots g_2 w_n g_2^{-1}\ldots g_{n-1}^{-1}$$ and  
$$
\begin{array}{rll}
b_n & = & g_{n-1}(g_{n-2}\ldots g_2w_{n-1}g_2^{-1}\ldots g_{n-2}^{-1}) \\
& &  g_{n-2}  (g_{n-3}\ldots g_2 w_{n-2} g_2^{-1}\ldots g_{n-3}^{-1}) g_{n-3}\ldots \\
& & g_4(g_3g_2w_4g_2^{-1}g_3^{-1}) g_3 (g_2w_3g_2^{-1})(g_{2}w_{2}g^{-1}_{2}) g_3^{-1}\ldots g_{n-1}^{-1}.
\end{array}
$$ 
Clearly, $a_{n} \in H_{n}$ and $b_{n} \in I_{n-1}$. Now, the initial equation can be written as $$g_na_nb_ng_n^{-1}=y$$ or $$(g_na_nb_ng_n^{-1}b_n^{-1}) b_n=y.$$ Since $H_{n}$ is a normal subgroup of $I_{n}$, we get $g_na_nb_ng_n^{-1}b_n^{-1}\in H_n$. Now, we may repeat the above decomposition as follows. Write $$a_{n-1}=g_{n-2}\ldots g_2 w_{n-1} g_2^{-1}\ldots g_{n-2}^{-1}\in H_{n-1}$$ and 
$$
\begin{array}{rll}
b_{n-1} & = & g_{n-2}(g_{n-3}\ldots g_2 w_{n-2}g_2^{-1}\ldots g_{n-3}^{-1}) g_{n-3} \\
& & \ldots g_4 (g_3g_2w_4g_2^{-1}g_3^{-1})g_3(g_2w_3g_2^{-1})(g_{2}w_{2}g^{-1}_{2})g_3^{-1}\ldots g_{n-2}^{-1}\in I_{n-2}.
\end{array}
$$
So, our equation has the form $$(g_na_nb_ng_n^{-1}b_n^{-1})(g_{n-1}a_{n-1}b_{n-1}g_{n-1}^{-1}b_{n-1}^{-1}) b_{n-1}=y,$$ with $g_na_nb_ng_n^{-1}b_n^{-1}\in H_n$ and $g_{n-1}a_{n-1}b_{n-1}g_{n-1}^{-1}b_{n-1}^{-1}\in H_{n-1}$. Repeating the above we will have, in finitely many steps,
$$(g_na_nb_ng_n^{-1}b_n^{-1}) (g_{n-1}a_{n-1}b_{n-1}g_{n-1}^{-1}b_{n-1}^{-1}) \ldots (g_3a_3b_3g_3^{-1}b_3^{-1})(g_2w_2g_2^{-1})=y,$$ where $g_ia_ib_ig_i^{-1}b_i^{-1}\in H_i$ for all $i\in\{3,\ldots,n\}$ and $g_2w_2g_2^{-1}\in H_2$. So, the word in the left is in normal form. Thus, for the two words to be equal we must have $g_2w_2g_2^{-1}=z_2$. But that  
is the conjugacy problem in $I_2 = {\rm Inn}(F_{2})\cong F_2$,  which is well known to be decidable (see, for example, \cite[Proposition 2.14]{lysc}). If such a $g_2$ doesn't exist, then $x,y$ are not conjugate and we are done. If such a $g_2$ exists, then we must also have $$g_3a_3b_3g_3^{-1}b_3^{-1}=z_3$$ or equivalently,
$$g_3  a_{3} (g_2w_2g_2^{-1}) g^{-1}_3 (g_2w_2^{-1}g_2^{-1})=z_3.$$ Since $g_2w_2g_2^{-1}\in H_2$ acts as an automorphism, say $\phi$, on $H_3$, we have $$g_3 a_{3} \phi(g_3^{-1})=z_3.$$ But this is the twisted conjugacy problem in $H_3 = {\rm Inn}(F_{3})\cong F_3$, which is decidable (see, for example, \cite{bmmv}). Hence, if such a $g_3$ does not exist, then $x,y$ are not conjugate. If it exists, we continue the same way. So, in a finite number of steps, either we fail to find some $g_i$ and so, $x,y$ are not conjugate, or we complete the algorithm by demonstrating the appropriate $g$. \qed

\vskip .120 in

In fact, the above proof shows the following more general result. 

\begin{cor}
Let $G$ be the iterated semi-direct product $G=(\ldots(G_n\rtimes G_{n-1})\rtimes \ldots \rtimes G_3) \rtimes G_2$ such that the group $G_{j}$ normalizes $G_{i}$ for all $i > j$, $i, j \in \{2, \ldots, n\}$. Then, the conjugacy problem is decidable in $G$ if the twisted conjugacy problem is decidable in $G_i$ for all $i=2,\ldots,n$.  
\end{cor}

\section{The Lie algebra of $I_n$}
Throughout this paper, by \lq \lq Lie algebra\rq \rq,  we mean a Lie
algebra over the ring of integers $\mathbb{Z}$ and we use the left-normed convention for group commutators and Lie commutators.  Let $G$ be a
group. For a positive integer $c$, let $\g_c(G)$ be the $c$-th term of
the lower central series of $G$. We point out that $\gamma_{2}(G) = G^{\prime}$; that is, the commutator subgroup of $G$. Write ${\rm gr}_{c}(G) = \g_c(G)/\g_{c+1}(G)$ for $c \geq 1$. The (restricted) direct sum of
the quotients ${\rm gr}_{c}(G)$ is the {\it associated graded
Lie algebra} of $G$, ${\rm gr}(G) = \bigoplus_{c \geq 1}
{\rm gr}_{c}(G)$. The Lie bracket multiplication in $\gr(G)$ is $
[a \gamma_{c+1}(G), b \gamma_{d+1}(G)] = [a, b]
\gamma_{c+d+1}(G)$, 
where $a \gamma_{c+1}(G)$ and $b \gamma_{d+1}(G)$ are the
images of the elements $a \in \gamma_{c}(G)$ and $b \in
\gamma_{d}(G)$ in the quotient groups ${\rm gr}_{c}(G)$ and
${\rm gr}_{d}(G)$, respectively, and $[a, b] \gamma_{c+d+1}(G)$
is the image of the group commutator $[a, b]$ in the quotient
group ${\rm gr}_{c+d}(G)$. Multiplication is then extended to
${\rm gr}(G)$ by linearity.  

In the present section, we investigate the Lie algebra ${\rm gr}(I_{n})$ of $I_n$ and we give a presentation for it. For $\mathbb{Z}$-submodules $A$ and $B$ of any Lie algebra, let $[A,B]$ be the $\mathbb{Z}$-submodule spanned by
$[a,b]$ where $a \in A$ and $b \in B$. Furthermore, $B \wr A$
denotes the $\mathbb{Z}$-submodule defined by $
B \wr A = B + [B,A] + [B,A,A] + \cdots$. Also, for any non-empty subsets ${\mathcal U}$ and ${\mathcal V}$ of any Lie algebra, we denote by $[{\mathcal V},{\mathcal U},\underbrace{{\mathcal U},\ldots, {\mathcal U}}_{a},\underbrace{{\mathcal V},\ldots, {\mathcal V}}_{b}]$ the set of (Lie) commutators where the first element belongs to ${\mathcal V}$, the second to ${\mathcal U}$ and then $a$ consecutive elements belong to ${\mathcal U}$ and the last $b$ consecutive elements belong to ${\mathcal V}$. For brevity,  this is also denoted by $[{\mathcal V},{\mathcal U},~_a{\mathcal U},~_b{\mathcal V}]$. For a free $\mathbb{Z}$-module $A$, let $L(A)$ be the free
Lie algebra on $A$, that is, the free Lie algebra on $\mathcal A$, where
$\mathcal A$ is an arbitrary $\mathbb{Z}$-basis of $A$. Thus, we may
write $L(A) = L({\mathcal A})$. For a positive integer $c$, let
$L^{c}(A)$ denote the $c$th homogeneous component of $L(A)$. It
is well-known that $
L(A) = \bigoplus_{c \geq 1}L^{c}(A)$. 

The following
result is a version of Lazard's "Elimination Theorem" (see
\cite[Chapter 2, Section 2.9, Proposition 10]{bour}). In the form
written here, it is a special case of \cite[Lemma 2.2]{bks1} (see, also, \cite[Section 2.2]{mp}). 

\begin{lemma}[Elimination Theorem]\label{le2}
Let $U$ and $V$ be free $\mathbb{Z}$-modules with $U\cap V=\{0\}$, and consider the
free Lie algebra $L(U \oplus V)$. Then, $U$ and $V \wr U$ freely
generate Lie subalgebras $L(U)$ and $L(V \wr U)$, and there is a
$\mathbb{Z}$-module decomposition $L(U \oplus V) = L(U) \oplus L(V
\wr U)$. Furthermore, $
V \wr U = V \oplus [V,U] \oplus [V, U, U] \oplus \cdots$ 
and, for each $n \geq 0$, there is a $\mathbb{Z}$-module isomorphism $
[V,~_nU] \cong V
\otimes \underbrace{U \otimes \cdots \otimes U}_{n}$ under which $[v,u_{1}, \ldots, u_{n}]$ corresponds to $v \otimes u_{1} \otimes \cdots \otimes u_{n}$ for all $v \in V$ and all $u_{1}, \ldots u_{n} \in U$. 
\end{lemma}

As a consequence of Lemma \ref{le2}, we have the following result.

\begin{corollary}\label{c1}
For free $\mathbb{Z}$-modules $U$ and $V$, we write $L(U \oplus V)$ for the free
Lie algebra on $U \oplus V$. Then, there is a
$\mathbb{Z}$-module decomposition $L(U \oplus V) = L(U) \oplus L(V) \oplus L(W)$,
where $W = W_{2} \oplus W_{3} \oplus \cdots$ such that, for all $m
\geq 2$, $W_{m}$ is the direct sum of submodules $[V,U, ~_{a}U,~_{b}V]$ with $a+b = m-2$ and $a, b \geq 0$. Each $[V,U,~_aU,~_bV]$ is isomorphic to $V \otimes U \otimes \underbrace{U \otimes \cdots \otimes U}_{a} \otimes \underbrace{V \otimes \cdots \otimes V}_{b}$ as
$\mathbb{Z}$-module. Furthermore, $L(W)$ is the ideal of $L(U \oplus V)$ generated by the module
$[V,U]$.
\end{corollary}

\subsection{A decomposition of the free Lie algebra}
Let $L = L({\mathcal Y})$ be the free Lie algebra on the finite set $\mathcal Y$ and decompose ${\mathcal Y}$ into a disjoint union ${\mathcal Y} = {\mathcal U} \cup {\mathcal V}$ of finite non-empty subsets ${\mathcal U} = \{u_{1}, \ldots, u_{\kappa}\}$ and ${\mathcal V} = \{v_{1}, \ldots, v_{\lambda}\}$. Also, define $u_{1} < \cdots < u_{\kappa} < v_{1} < \cdots < v_{\lambda}$. Let $U$ and $V$ denote the free ${\mathbb{Z}}$-modules with bases ${\mathcal U}$ and ${\mathcal V}$, respectively. Since the $\mathbb{Z}$-module generated by ${\mathcal Y}$ is equal to $U \oplus V$ and $L({\mathcal Y})$ is free on the set ${\mathcal Y}$, $L$ is free on the set ${\mathcal U} \cup {\mathcal V}$ and we have $L=L(U\oplus V)$.  By Lemma \ref{le2} and Corollary \ref{c1}, we have 
$$
\begin{array}{lll}
L & = & L(U \oplus V) \\
& = & L(U) \oplus L(V \wr U) \\
& = & L(U) \oplus L(V) \oplus L(W),
\end{array}
$$
where $W = W_{2} \oplus W_{3} \oplus \cdots$ such that, for all $m \geq 2$, $W_{m} = \bigoplus_{a+b=m-2}[V,U,~_{a}U, ~_{b}V]$. 
Furthermore, $L(V \wr U)$  and $L(W)$ are the ideals in $L$ generated by the modules $V$ and $[V,U]$, respectively. In particular, $L(W)$ is the ideal in $L$ generated by the natural $\mathbb{Z}$-basis $
[{\mathcal V}, {\mathcal U}] = \{[v_{i}, u_{j}]: 1 \leq i \leq \lambda, 1 \leq j \leq \kappa\}$  
of $[V,U]$. Let ${\mathcal X}_{V,U}$ be the natural $\mathbb{Z}$-basis of $V \wr U$. That is, 
$$
{\mathcal X}_{V,U} = {\mathcal V} \cup( \bigcup_{a \geq 1}[{\mathcal V},~_{a}{\mathcal U}]),
$$ 
where $[{\mathcal V},~_{a}{\mathcal U}] = \{[v, z_{1}, \ldots, z_{a}]: v \in {\mathcal V}, z_{1}, \ldots, z_{a} \in {\mathcal U}\}$  
is the natural $\mathbb{Z}$-basis of the module $[V,~_{a}U]$. We point out that the set 
$$
{\mathcal W} = [{\mathcal V}, {\mathcal U}] \cup (\bigcup_{a+b \geq 1 \atop a, b \geq 0}[{\mathcal V}, {\mathcal U},~_{a}{\mathcal U},~_{b}{\mathcal V}]),
$$
where $[{\mathcal V}, {\mathcal U},~_{a}{\mathcal U},~_{b}{\mathcal V}] = \{[v,z,z_{1}, \ldots, z_{a}, v_{1}, \ldots,v_{b}]: v, v_{1}, \ldots, v_{b} \in {\mathcal V}, z, z_{1}, \ldots, z_{a} \in {\mathcal U}\}$, is the natural $\mathbb{Z}$-basis of $W$.

\begin{lemma}\label{le3}
With the above notation, let $\phi$ be any automorphism of the free $\mathbb{Z}$-module $[V,U]$ and let $\psi_{2}: [{\mathcal V}, {\mathcal U}] \rightarrow L(V \wr U)$ be the map given by $\psi_{2}([v_{i},u_{j}]) = \phi([v_{i},u_{j}]) + w_{i,j}$, 
where $w_{i,j} \in L^{2}(V)$ for all $1 \leq i \leq \lambda$ and $1 \leq j \leq \kappa$. Furthermore, for every $a \geq 3$, let $\psi_{a}$ be the mapping from $[{\mathcal V},{\mathcal U}, ~_{(a-2)}{\mathcal U}]$ into $L(V \wr U)$ satisfying the conditions $
\psi_{a}([v,z,z_{1}, \ldots, z_{a-2}]) = [\psi_{2}([v,z]),z_{1}, \ldots, z_{a-2}]$  
for all $v \in {\mathcal V}$ and $z,z_{1}, \ldots, z_{a-2} \in {\mathcal U}$. Define $\Psi: {\mathcal X}_{V,U} \rightarrow L(V \wr U)$ to be the induced map with
$\Psi(v) = v$ for all $v \in {\mathcal V}$ and, for $a \geq 2$, $\Psi(v) = \psi_{a}(v)$ for all $v \in [{\mathcal V}, {\mathcal U}, ~_{(a-2)}{\mathcal U}]$. Then, 
\begin{enumerate}

\item $\psi_{2}$ extends to a $\mathbb{Z}$-linear monomorphism from $[V,U]$ into $L(V\wr U)$. 

\item For every $a \geq 3$, there exists an automorphism $\phi_{a}$ of the free $\mathbb{Z}$-module $[V,U, ~_{(a-2)}U]$ such that, for all $v \in {\mathcal V}$ and $z,z_{1}, \ldots, z_{a-2} \in {\mathcal U}$,  
$$
\psi_{a}([v,z,z_{1}, \ldots, z_{a-2}]) = \phi_{a}([v,z,z_{1}, \ldots, z_{a-2}]) + w_{a},
$$
where $w_{a} \in L({\mathcal V} \cup (\bigcup_{1 \leq \beta \leq a-2}[{\mathcal V}, ~_{\beta}{\mathcal U}]))$. 

\item $\Psi$ extends to a (Lie algebra) automorphism of $L(V \wr U)$.  

\end{enumerate} 
\end{lemma}

\pf \begin{enumerate}

\item Since $[{\mathcal V}, {\mathcal U}]$ is a $\mathbb{Z}$-basis of $[V,U]$, $\psi_{2}$ extends to a $\mathbb{Z}$-linear mapping, denoted $\psi_{2}$ again, from $[V,U]$ into $L(V \wr U)$. We point out that $\psi_{2}([V,U]) \subseteq L^{2}(V) \oplus [V,U]$. Since $L^{2}(V) \cap [V,U] = \{0\}$, $\phi$ is an automorphism of $[V,U]$ and the set $[{\mathcal V},{\mathcal U}]$ is a $\mathbb{Z}$-basis of $[V,U]$, we have $\psi_{2}$ is a $\mathbb{Z}$-linear monomorphism of $[V,U]$ into $L(V \wr U)$.

\item Let $a \geq 3$. Since $[{\mathcal V}, {\mathcal U}, ~_{(a-2)}{\mathcal U}]$ is a $\mathbb{Z}$-basis for the $\mathbb{Z}$-module $[V,U, ~_{(a-2)}U]$ and $\psi_{2}$ is a $\mathbb{Z}$-linear mapping, we have $\psi_{a}$ extends to a $\mathbb{Z}$-linear mapping from $[V,U, ~_{(a-2)}U]$ into $L(V \wr U)$. For $a \geq 3$, let $\phi_{a}$ be the mapping from $[{\mathcal V},{\mathcal U}, ~_{(a-2)}{\mathcal U}]$ into $[V,U, ~_{(a-2)}U]$ satisfying the conditions $\phi_{a}([v,z,z_{1}, \ldots,$  $z_{a-2}]) = [\phi([v,z]), z_{1}, \ldots, z_{a-2}]$ 
for all $v \in {\mathcal V}$ and $z, z_{1}, z_{a-2} \in {\mathcal U}$. Since $\phi$ is a $\mathbb{Z}$-linear mapping, $\phi_{a}$ extends to a $\mathbb{Z}$-linear mapping of $[V,U, ~_{(a-2)}U]$. Let $\chi_{a}$ be the mapping from $[{\mathcal V},{\mathcal U}, ~_{(a-2)}{\mathcal U}]$ into $[V,U, ~_{(a-2)}U]$ satisfying the conditions  $\chi_{a}([v,z,z_{1}, \ldots, z_{a-2}]) = [\phi^{-1}([v,z]), z_{1}, \ldots, z_{a-2}]$ 
for all $v \in {\mathcal V}$, $z,z_{1}, \ldots, z_{a-2} \in {\mathcal U}$, where $\phi^{-1}$ is the inverse of $\phi$. We point out that the elements $z_{1}, \ldots, z_{a-2} \in {\mathcal U}$ are invariant under $\phi_{a}$ and $\chi_{a}$ in the above expressions. Since $\phi^{-1}$ is the inverse of $\phi$, $z_{1}, \ldots, z_{a-2}$ are invariant under $\phi_{a}$ and $\chi_{a}$ and $[{\mathcal V}, {\mathcal U}, ~_{(a-2)}{\mathcal U}]$ is a $\mathbb{Z}$-basis of $[V,U, ~_{(a-2)}U]$, we have $\chi_{a}$ is the inverse of $\phi_{a}$. Therefore, $\phi_{a}$ is an automorphism of $[V,U, ~_{(a-2)}U]$. Now, for $v \in {\mathcal V}$, $z,z_{1}, \ldots, z_{a-2} \in {\mathcal U}$,  
$$
\begin{array}{ccl}
\psi_{a}([v,z,z_{1}, \ldots, z_{a-2}]) & = &[\psi_{2}([v,z]),z_{1}, \ldots, z_{a-2}] \\
& = & [\phi([v,z]),z_{1}, \ldots, z_{a-2}] + [w_{v,z}, z_{1}, \ldots, z_{a-2}] \\
& = & \phi_{a}([v,z,z_{1}, \ldots, z_{a-2}]) + [w_{v,z}, z_{1}, \ldots, z_{a-2}], 
\end{array}
$$
where $w_{v,z} \in L^{2}(V)$. By using the Jacobi identity in the form $[x,y,z] = [x,z,y] + [x, [y,z]]$, 
we may show that $$[w_{v,z}, z_{1}, \ldots, z_{a-2}] \in L({\mathcal V} \cup (\bigcup_{1 \leq \beta \leq a-2}[{\mathcal V}, ~_{\beta}{\mathcal U}])).$$ Therefore,
for all $v \in {\mathcal V}$ and $z,z_{1}, \ldots, z_{a-2} \in {\mathcal U}$, 
$$
\psi_{a}([v,z,z_{1}, \ldots, z_{a-2}]) = \phi_{a}([v,z,z_{1}, \ldots, z_{a-2}]) + w_{a},
$$ 
where $w_{a} \in L({\mathcal V} \cup (\bigcup_{1 \leq \beta \leq a-2}[{\mathcal V}, ~_{\beta}{\mathcal U}]))$.

\item Since ${\mathcal X}_{V,U}$ is a free generating set of $L(V \wr U)$, $\Psi$ extends to an endomorphism of $L(V \wr U)$. Since $\Psi(v) = v$ for all $v \in {\mathcal V}$, we obtain from Lemma \ref{le3}~(1)-(2) and \cite[Lemma 2.1]{bks2} that $\Psi$ is an automorphism of $L(V \wr U)$. \qed

\end{enumerate} 
Since $L(W)$ is a free Lie subalgebra of $L(V \wr U)$ on $W$ and $\Psi$ is an automorphism of $L(V \wr U)$, we have $\Psi(L(W))$ is a free Lie subalgebra of $L(V \wr U)$. In fact, $\Psi(L(W)) = L(\Psi(W))$,   
that is, $\Psi(L(W))$ is a free Lie algebra on $\Psi(W)$. Furthermore, since $\Psi$ is an automorphism of $L(V \wr U)$ and $L(W)$ is an ideal of $L(V \wr U)$, we obtain $\Psi(L(W)) = L(\Psi(W))$ is an ideal in $L(V \wr U)$. We point out that 
$$
\begin{array}{rll}
L(V \wr U) & = & \Psi(L(V \wr U)) \\
({\rm By~ Corollary}~ \ref{c1}) & = & \Psi(L(V) \oplus L(W)) \\
(\Psi~{\rm automorphism}) & = & \Psi(L(V)) \oplus \Psi(L(W)) \\
& = & L(V) \oplus L(\Psi(W))
\end{array}
$$
and so, $
\begin{array}{lll}
L & = & L(U) \oplus L(V) \oplus L(\Psi(W)).
\end{array}$
Write 
$$
{\mathcal C} = \Psi({\mathcal W}) = \Psi([{\mathcal V}, {\mathcal U}]) \cup (\bigcup_{a+b \geq 1 \atop a, b \geq 0}[\Psi([{\mathcal V}, {\mathcal U}]), ~_{a}{\mathcal U}, ~_{b}{\mathcal V}]). 
$$
Since $\mathcal W$ is a $\mathbb{Z}$-basis for $W$ and $\Psi$ is an automorphism of $L(V \wr U)$, we have $\mathcal C$ is a $\mathbb{Z}$-basis of $\Psi(W)$ and so, $L(\Psi(W)) = L({\mathcal C})$. For a positive integer $m$, with $m \geq 2$, let 
$$
W_{m,\Psi} = \bigoplus_{a, b \geq 0 \atop a+b = m-2}[\Psi([V,U]), ~_{a}U, ~_{b}V].
$$
In particular, for $m \geq 2$,  the set 
$$
{\mathcal W}_{m,\Psi} = \bigcup_{a, b \geq 0 \atop a+b = m-2}[\Psi([{\mathcal V},{\mathcal U}]), ~_{a}{\mathcal U}, ~_{b}{\mathcal V}]
$$
is a natural $\mathbb{Z}$-basis for $W_{m,\Psi}$ and so, $\Psi(W) = \bigoplus_{m \geq 2} W_{m,\Psi}$. Obviously,  ${\mathcal C} = \bigcup_{m \geq 2}{\mathcal W}_{m,\Psi}$.   
Furthermore, for $m \geq 2$, we write $
L^{m}_{\rm grad}(\Psi(W))$ for $L^{m} \cap L(\Psi(W))$. That is, $L^{m}_{\rm grad}(\Psi(W))$ is the $\mathbb{Z}$-submodule of $L^{m}$ spanned by all Lie commutators of the form $[w_{1}, \ldots, w_{\kappa}]$, with $\kappa \geq 1$, $w_{i} \in {\mathcal W}_{m(i),\Psi}$ and $m(1) + \cdots + m(\kappa) = m$. It is easily checked that 
$$
L(\Psi(W)) = \bigoplus_{m \geq 2} L^{m}_{\rm grad}(\Psi(W)).
$$   

\begin{proposition}\label{pr1}
With the above notation, $L(\Psi(W))$ is an ideal in $L$.
\end{proposition}

\pf To show that the Lie subalgebra $L(\Psi(W))$ of $L$ is an ideal in $L$, it is enough to show that $[w,u] \in L(\Psi(W))$ for all $w \in L(\Psi(W))$ and $u \in L$. Since $L(\Psi(W))$ is an ideal in $L(V \wr U)$ and by Lemma \ref{le2}, it is enough to show that $[w,u] \in L(\Psi(W))$ for all $w \in L(\Psi(W))$ and $u \in L(U)$. Since any $u \in L(U)$ is written as a $\mathbb{Z}$-linear combination of Lie commutators $[u_{i_{1}}, \ldots, u_{i_{\kappa}}]$, with $\kappa \geq 1$ and $u_{i_{1}}, \ldots, u_{i_{\kappa}} \in {\mathcal U}$, and by using the Jacobi identity in the form $[x,[y,z]] = [x,y,z] - [x,z,y]$, it is enough to show that $[w, u_{i_{1}}, \ldots, u_{i_{k}}] \in L(\Psi(W))$ for all $w \in L(\Psi(W))$ and $u_{i_{1}}, \ldots, u_{i_{k}} \in {\mathcal U}$. Since $\mathcal C$ is a (free) generating set of $L(\Psi(W)) = L({\mathcal C})$, any $w \in L(\Psi(W))$ is written as a $\mathbb{Z}$-linear combination of Lie commutators $[w_{1}, \ldots, w_{\lambda}]$ with $w_{1}, \ldots, w_{\lambda} \in {\mathcal C}$. By the linearity of the Lie bracket, it is enough to show that $[w_{1}, \ldots, w_{\lambda}, u_{i_{1}}, \ldots, u_{i_{\kappa}}] \in L(\Psi(W))$ for all $w_{1}, \ldots, w_{\lambda} \in {\mathcal C}$ and $u_{i_{1}}, \ldots, u_{i_{\kappa}} \in {\mathcal U}$. By using induction on $\lambda$, the Jacobi identity in the form $[x,y,z] = [x,z,y] + [x,[y,z]]$ and the fact that $L(\Psi(W))$ is a Lie algebra, it suffices to show that $[w, u_{i_{1}}, \ldots, u_{i_{\kappa}}] \in L(\Psi(W))$ for all $w \in {\mathcal C}$ and $u_{i_{1}}, \ldots, u_{i_{\kappa}} \in {\mathcal U}$. 

Notice that if  $w \in \Psi([{\mathcal V}, {\mathcal U}])$, then $[w, u_{i_{1}}, \ldots, u_{i_{\kappa}}] \in {\mathcal C}$ for all $u_{i_{1}}, \ldots, u_{i_{\kappa}} \in {\mathcal U}$. Thus, we assume that $w \in [\Psi([{\mathcal V}, {\mathcal U}]), ~_{a}{\mathcal U}, ~_{b}{\mathcal V}]$ for some integers $a, b \geq 0$  with $a + b \geq 1$ and $b \geq 1$. For elements $u, v_{1}, \ldots, v_{m} \in L$, we denote    
$$
[u, \underbrace{v_{1},\ldots, v_1}_{\mu_{1}}, \ldots, \underbrace{v_{m},\ldots,v_{m}}_{\mu_{m}}]=[u, ~_{\mu_{1}}v_1, \ldots, ~_{\mu_{m}}v_{m}].
$$ 
Rewrite $w = [u, v_{i_{1}}, \ldots, v_{i_{b}}]$, where $u \in [\Psi([{\mathcal V},{\mathcal U}]),~_{a}{\mathcal U}]$, $v_{i_{1}}, \ldots, v_{i_{b}} \in {\mathcal V}$, $a, b$ non-negative integers with $a + b \geq 1$ and $b \geq 1$. By  applying repeatedly the Jacobi identity in the form $[x,y,z] = [x,z,y] + [x, [y,z]]$, we get $w = [u, ~_{\mu_{1}}v_1, \ldots, ~_{\mu_{\lambda}}v_{\lambda}] + \tilde{w}$, 
where $\mu_{1} + \cdots + \mu_{\lambda} = b$ and  $\tilde{w} \in L(V \wr U)$. In particular, $\tilde{w}$ is a summand of elements of the form $[u, ~_{\mu^{\prime}_{1}}v_1, \ldots, ~_{\mu^{\prime}_{\beta}}v_{\beta}, z_{i_{1}}, \ldots, z_{i_{\gamma}}]$,   
where $\mu^{\prime}_{1}, \ldots, \mu^{\prime}_{\beta} \geq 0$, $\mu^{\prime}_{1} + \cdots +  \mu^{\prime}_{\beta} \leq b$ and $z_{i_{1}}, \ldots, z_{i_{\gamma}} \in L^{\prime}(V)$. Since $L(\Psi(W))$ is an ideal in $L(V \wr U)$, it suffices to consider the case where  $w = [u, ~_{\mu_{1}}v_1, \ldots, ~_{\mu_{\la}}v_{\la}]$, with $u \in [\Psi([{\mathcal V},{\mathcal U}]),~_{a}{\mathcal U}]$, $\mu_{1}, \ldots, \mu_{\la} \geq 0$, $\mu_{1} + \cdots + \mu_{\la} = b \geq 1$ and $v_{1}, \ldots, v_{\la} \in {\mathcal V}$.  Since $\phi$ is an automorphism of $[V,U]$, we may write, for $i \in \{1, \ldots, \la\}$, $j \in \{1, \ldots, \ka\}$, 
$$
\begin{array}{rllrr}
[v_{i}, u_{j}] & = & \sum_{r=1}^{\la} \sum_{s=1}^{\ka} b_{i,j,r,s} ~\phi([v_{r},u_{s}]) & & \\
& = & \sum_{r=1}^{\la} \sum_{s=1}^{\ka} b_{i,j,r,s} (\psi_{2}([v_{r},u_{s}]) + w^{\prime}_{r,s}) & & \\
& = & (\sum_{r=1}^{\la} \sum_{s=1}^{\ka} b_{i,j,r,s} \psi_{2}([v_{r},u_{s}])) + \tilde{w}_{i,j} & & \ \ \ \ \ \ \ (1)\\
\end{array}
$$
where $b_{i,j,r,s} \in \mathbb{Z}$ and  $\tilde{w}_{i,j} \in L^{2}(V)$. By using the Jacobi identity in the expression $[u, ~_{\mu_{1}}v_1, \ldots, ~_{\mu_{m}}v_{m}, u_{i_{1}}, \ldots, u_{i_{k}}]$, where $u \in [\Psi([{\mathcal V},{\mathcal U}]),~_{a}{\mathcal U}]$, $\mu_{1}, \ldots, \mu_{m} \geq 0$, $\mu_{1} + \cdots + \mu_{m} = b \geq 1$, $v_{1}, \ldots, v_{m} \in {\mathcal V}$ and $u_{i_{1}}, \ldots, u_{i_{k}} \in {\mathcal U}$, replacing each $[v_{i}, u_{j}]$ by the equation $(1)$ as many times as it is needed and since $L(\Psi(W))$ is an ideal in $L(V \wr U)$, we may show that $[u, ~_{\mu_{1}}v_1, \ldots, ~_{\mu_{m}}v_{m}, u_{i_{1}}, \ldots, u_{i_{k}}] \in L(\Psi(W))$ (see \cite[Example 1]{mps} for an example of the above procedure). Therefore, $L(\Psi(W))$ is an ideal in $L$. \qed

\vskip .120 in

Since 
$$
\begin{array}{rll}
L  & = & L(U) \oplus L(V) \oplus L(\Psi(W)) \\
& = & L(U) \oplus L(V) \oplus (\bigoplus_{m \geq 2} L^{m}_{\rm grad}(\Psi(W))),
\end{array}
$$
we have the following result.

\begin{corollary}\label{co2}
Let $\Psi$ be the Lie algebra automorphism of $L(V \wr U)$ defined above. Then, for every positive integer $m$, with $m \geq 2$, $
L^{m} = L^{m}(U) \oplus L^{m}(V) \oplus L^{m}_{\rm grad}(\Psi(W))$. 
\end{corollary}

\subsection{A further decomposition}

Let $L({\mathcal Y})$ be  the free Lie algebra on the finite set $\mathcal Y$ and decompose $\mathcal Y$ into a disjoint union ${\mathcal Y} = {\mathcal Y}_{2} \cup \cdots \cup {\mathcal Y}_{n}$, with $n \geq 3$, of finite non-empty subsets of the form ${\mathcal Y}_{i} = \{y_{i1}, y_{i2}, \ldots, y_{ik_{i}}\}$,  
$2 \leq i \leq n$ and $k_i\ge 1$. Moreover, define $y_{i1} < \cdots < y_{ik_{i}}$ and $y_{ik_{i}} < y_{(i+1)1}$ for all $i \in \{2, \ldots, n-1\}$. 
For $m \in \{2, \ldots, n\}$, let $Y_{m}$ be the free $\mathbb{Z}$-module with $\mathbb{Z}$-basis ${\mathcal Y}_{m}$. Thus, 
$$
L = L({\mathcal Y}) = L(Y_{2} \oplus \cdots \oplus Y_{n}).
$$
Then, 
$$
L^{2}({\mathcal Y}) = \left(\bigoplus_{i=2}^{n}L^{2}(Y_{i})\right) \oplus \left(\bigoplus_{2 \leq \kappa < \rho \leq n}[Y_{\rho}, Y_{\kappa}]\right).
$$
Write $W_{(2, \ldots, n),2} = \bigoplus_{2 \leq \kappa < \rho \leq n}[Y_{\rho}, Y_{\kappa}]$. For $r \in \{2, \ldots, n-1\}$, let $U_{r} = Y_{r} \oplus \cdots \oplus Y_{n}$. We point out that ${\mathcal{U}}_{r} = \bigcup_{\kappa = r}^{n}{\mathcal Y}_{\kappa}$ is a natural $\mathbb{Z}$-basis of $U_{r}$. It is easily verified that 
$$
W_{(2, \ldots, n),2} = \bigoplus_{r=2}^{n-1}[U_{r+1},Y_{r}].
$$
For $r \in \{2, \ldots, n-1\}$, $[U_{r+1},Y_{r}]$ is generated by the following standard free generating set $
[{\mathcal U}_{r+1}, {\mathcal Y}_{r}] = \{[y_{\mu \nu},y_{rl}]: r+1 \leq \mu \leq n, \nu = 1, \ldots, k_{\mu}, l = 1, \ldots, k_{r}\}$. Furthermore, $W_{(2, \ldots, n),2}$ is generated by the set $K = \bigcup_{r=2}^{n-1}[{\mathcal U}_{r+1}, {\mathcal Y}_{r}]$, which is a natural $\mathbb{Z}$-basis of $W_{(2, \ldots, n),2}$. Let $\phi_{r}$ be any automorphism of the free $\mathbb{Z}$-module $[U_{r+1},Y_{r}]$ and let $\psi_{2,r}: [{\mathcal U}_{r+1}, {\mathcal Y}_{r}] \rightarrow [U_{r+1},Y_{r}] \oplus L^{2}(U_{r+1})$ be the map given by $\psi_{2,r}([y_{\mu \nu},y_{rl}]) = \phi_{r}([y_{\mu \nu},y_{rl}]) + w_{\mu,\nu,r,l}$, where $w_{\mu,\nu,r,l} \in L^{2}(U_{r+1})$, $\mu = r+1, \ldots, n$, $\nu = 1, \ldots, k_{\mu}$ and $l = 1, \ldots, r$. We point out that 
$$
L^{2}(U_{r+1}) = \left(\bigoplus_{j=r+1}^{n}L^{2}(Y_{j})\right) \oplus \left(\bigoplus_{r+1 \leq \kappa < \rho \leq n}[Y_{\rho}, Y_{\kappa}]\right).
$$
Since $\phi_{r}$ is an automorphism and $[U_{r+1},Y_{r}] \cap L^{2}(U_{r+1}) = \{0\}$, we get $\psi_{2,r}$ is $1-1$. For $a \geq 3$, let $\psi_{a,r}$ be the mapping from $[{\mathcal U}_{r+1}, {\mathcal Y}_{r}, ~_{(a-2)}{\mathcal Y}_{r}]$ into $L(U_{r+1} \wr Y_{r})$ satisfying the conditions 
$$
\psi_{a,r}([u,z,z_{1}, \ldots, z_{a-2}]) = [\psi_{2,r}([u,z]), z_{1}, \ldots, z_{a-2}]
$$
for all $u \in {\mathcal U}_{r+1}$ and $z, z_{1}, \ldots, z_{a-2} \in {\mathcal Y}_{r}$. Let ${\mathcal X}_{U_{r+1},Y_{r}}$ be the natural $\mathbb{Z}$-basis of $U_{r+1} \wr Y_{r}$. That is, 
$$
{\mathcal X}_{U_{r+1},Y_{r}} = {\mathcal U}_{r+1} \cup( \bigcup_{a \geq 1}[{\mathcal U}_{r+1},~_{a}{\mathcal Y}_{r}]),
$$
where $
[{\mathcal U}_{r+1},~_{a}{\mathcal Y}_{r}] = \{[u, z_{1}, \ldots, z_{a}]: u \in {\mathcal U}_{r+1}, z_{1}, \ldots, z_{a} \in {\mathcal Y}_{r}\}$  
is the natural $\mathbb{Z}$-basis of the module $[U_{r+1},~_{a}Y_{r}]$. Furthermore, we define $\Psi_{r}: {\mathcal X}_{U_{r+1},Y_{r}} \rightarrow L(U_{r+1} \wr Y_{r})$ to be the map with
$\Psi_{r}(u) = u$ for all $u \in {\mathcal U}_{r+1}$ and, for $a \geq 2$, $\Psi_{r}(u) = \psi_{a,r}(u)$ for all $u \in [{\mathcal U}_{r+1}, {\mathcal Y}_{r}, ~_{(a-2)}{\mathcal Y}_{r}]$. By Lemma \ref{le3}~(3), $\Psi_{r}$ extends to a (Lie algebra) automorphism of $L(U_{r+1} \wr Y_{r})$. By applying Lemma \ref{le2} on  $L(U_{r})$, we have
$$
\begin{array}{lll}
L(U_{r}) & = & L(Y_{r}) \oplus L(U_{r+1} \wr Y_{r}) \\
& = & L(Y_{r}) \oplus L(U_{r+1}) \oplus L(W_{(r)}),
\end{array}
$$
where $W_{(r)} = W_{(r),2} \oplus W_{(r),3} \oplus \cdots$ such that, for all $\kappa \geq 2$, 
$$
W_{(r),\kappa} = \bigoplus_{a,b\ge 0 \atop a+b=\kappa-2}[U_{r+1},Y_{r},~_{a}Y_{r}, ~_{b}U_{r+1}].
$$
Furthermore, $L(U_{r+1} \wr Y_{r})$  and $L(W_{(r)})$ are the ideals in $L(U_{r})$ generated by the modules $U_{r+1}$ and $[U_{r+1},Y_{r}]$, respectively. By Proposition \ref{pr1}, $
T_{r} = L(\Psi_{r}(W_{(r)}))$  
is an ideal of $L(U_{r})$. Thus, for $r \in \{2, \ldots, n-1\}$, 
$$
[T_{r}, L(U_{r})] \subseteq T_{r}. \eqno(2)
$$
Furthermore, for $r \in \{2, \ldots, n-1\}$, let
$$
{\mathcal C}_{(r)} = \bigcup_{\kappa \geq 2}{\mathcal C}_{(r),\kappa},
$$
where, for $\kappa \geq 2$,   
$$
{\mathcal C}_{(r),\kappa} = \bigcup_{a, b \geq 0 \atop a+b=\kappa -2}[\psi_{2,r}([{\mathcal U}_{r+1}, {\mathcal Y}_{r}]), ~_{a}{\mathcal Y}_{r}, ~_{b}{\mathcal U}_{r+1}] 
$$
and ${\mathcal U}_{r+1} = {\mathcal Y}_{r+1} \cup \cdots \cup {\mathcal Y}_{n}$. 
Thus, for $r \in \{2, \ldots, n-1\}$, $
T_{r} = L(\Psi_{r}(W_{(r)})) = L({\mathcal C}_{(r)})$ and  $
L(U_{r}) = L(Y_{r}) \oplus L(U_{r+1}) \oplus T_{r}$. 
We define, for $r \in \{2, \ldots, n-1\}$, 
$$
\tilde{T}_{r} = \bigoplus_{s=2}^{r}T_{s}.
$$
We point out that, for $r \geq 3$, $
\tilde{T}_{r} = T_{r} \oplus \tilde{T}_{r-1}$. 

\begin{proposition}\label{pr2}
With the above notation, for every $r \in \{3, \ldots, n-1\}$,
\begin{enumerate}
\item $
L = \left(\bigoplus_{j=2}^{r-1}L(Y_{j})\right) \oplus L(U_{r}) \oplus \tilde{T}_{r-1}$. 
In particular, for $r = n-1$, we have  
$L = \left(\bigoplus_{i=2}^{n}L(Y_{i})\right) \oplus \tilde{T}_{n-1}$. 
\item $
[T_{r},L] \subseteq \left(\sum_{j=2}^{r-1}[T_{r},L(Y_{j})]\right) + T_{r} + [T_{r},\tilde{T}_{r-1}]$.  
\end{enumerate}

\end{proposition}

\pf 
(1)  We induct on $r$. Let $r = 3$. By Lemma \ref{le2}, $
L = L(Y_{2}) \oplus L(U_{3} \wr Y_{2})$. 
By Lemma \ref{le3} (3), $\Psi_{2}$ extends to a Lie automorphism of $L(U_{3} \wr Y_{2})$ and so, 
$$
\begin{array}{rll}
L & = & L(Y_{2}) \oplus L(U_{3} \wr Y_{2}) \\
& = & L(Y_{2}) \oplus \Psi_{2}(L(U_{3} \wr Y_{2})) \\
& = & L(Y_{2}) \oplus \Psi_{2}(L(U_{3})) \oplus L(W_{(2)})) \\
& = & L(Y_{2}) \oplus \Psi_{2}(L(U_{3})) \oplus \Psi_{2}(L(W_{(2)})) \\
& = & L(Y_{2}) \oplus L(U_{3}) \oplus \Psi_{2}(L(W_{(2)})) \\
& = & L(Y_{2}) \oplus L(U_{3}) \oplus T_{2}. 
\end{array}
$$
We assume that our claim is valid for $\kappa - 1$, with $4 \leq \kappa \leq n-1$. Thus, 
$$
L = \left(\bigoplus_{j=2}^{\kappa-2}L(Y_{j})\right) \oplus L(U_{\kappa - 1}) \oplus \left(\bigoplus_{j=2}^{\kappa-2}T_{j}\right). 
$$
By Lemma \ref{le2}, applied on $L(U_{\kappa-1})$, we have 
$
L(U_{\kappa-1}) = L(Y_{\kappa-1}) \oplus L(U_{\kappa} \wr Y_{\kappa-1}).
$
As before, 
by Lemma \ref{le3} (3), $\Psi_{\kappa-1}$ extends to a Lie automorphism of $L(U_{\kappa} \wr Y_{\kappa-1})$ and so, by Proposition \ref{pr1}, 
$
T_{\kappa-1} = L(\Psi_{\kappa-1}(W_{(\kappa-1)}))
$ 
is an ideal of $L(U_{\kappa-1})$. But   
$$
L(U_{\kappa-1}) = L(Y_{\kappa-1}) \oplus L(U_{\kappa}) \oplus T_{\kappa-1}. 
$$
Thus, we have the following decomposition of $L$,
$$
L = \left(\bigoplus_{j=2}^{\kappa-1}L(Y_{j})\right) \oplus L(U_{\kappa}) \oplus \left(\bigoplus_{j=2}^{\kappa-1}T_{j}\right) 
$$
and so, we obtain the required result.

\bigskip

\noindent (2) By Proposition \ref{pr2}~(1), for $r \in \{3, \ldots, n-1\}$,   
$$
L = \left(\bigoplus_{j=2}^{r-1}L(Y_{j})\right) \oplus L(U_{r}) \oplus \tilde{T}_{r-1}. \eqno(3) 
$$
For any $r \in \{3, \ldots, n-1\}$ and by using the equations (2) and (3), we get the required result. \qed 

By Corollary \ref{co2}, for any $r \in \{3, \ldots, n-1\}$ and  $m \geq 2$, 
$$
L^{m}(U_{r}) = L^{m}(Y_{r}) \oplus L^{m}(U_{r+1}) \oplus L^{m}_{\rm grad}(\Psi_{r}(W_{(r)})). 
$$
Since $\tilde{T}_{n-1} = \bigoplus_{s=2}^{n-1}L_{\rm grad}(\Psi_{s}(W_{(s)}))$, we have, for any $m \geq 2$, 
$$
\tilde{T}_{n-1}^{m} =  \bigoplus_{s=2}^{n-1}L^{m}_{\rm grad}(\Psi_{s}(W_{(s)})).
$$  
By Proposition \ref{pr2} (1) and the above equation, we have the following result.

\begin{corollary}\label{co4}
For any positive integer $m$, with $m \geq 2$, $
L^{m} = \left(\bigoplus_{i=2}^{n}L^{m}(Y_{i})\right) \oplus \tilde{T}_{n-1}^{m}$. 
\end{corollary} 

Let $J$ be the ideal of $L$ 
generated by the set ${\mathcal J}$ of all elements of the form $\psi_{2,r}
([y_{\mu \nu},y_{rl}])$, $r = 2, \ldots, n-1$, $\mu = r+1, \ldots, n-1$, $\nu = 1, \ldots, k_{\mu}$ and $l = 1, \ldots, r$. Since $J$ is a homogeneous ideal of $L$, we have $J = \bigoplus_{m \geq 2} J^{m}$, where $J^{m} = J \cap L^{m}$. 

\begin{lemma}\label{l3}
With the above notation, if $\tilde{T}_{n-1}$ is an ideal of $L$, then   
$\tilde{T}_{n-1} = J$ and so, $L=(\bigoplus_{i=2}^nL(Y_i))\oplus J$. In particular, for any positive integer $m$, with $m \geq 2$, $
J^{m} = \tilde{T}_{n-1}^{m} =  \bigoplus_{s=2}^{n-1}L^{m}_{\rm grad}(\Psi_{s}(W_{(s)}))$ 
and $
L^{m} = \left(\bigoplus_{i=2}^{n}L^{m}(Y_{i})\right) \oplus J^{m}$.  
\end{lemma}

\pf Let $\tilde{T}_{n-1}$ be an ideal of $L$. We point out that ${\mathcal{J}} \subseteq \tilde{T}_{n-1} \subseteq J$. Since $\tilde{T}_{n-1}$ is an ideal in $L$, we obtain $J \subseteq \tilde{T}_{n-1}$, and so, $J = \tilde{T}_{n-1}$. Hence, by Corollary \ref{co4} and since $J$ is homogeneous, we obtain the required result. \qed

\subsection{The case of $I_n$}\label{3.3}

Let us now consider the case of the subgroup $I_n$ of $M_{n}$, with $n \geq 3$. Let $L({\mathcal Y})$ be the free Lie algebra on the set ${\mathcal Y}$ which is decomposed as ${\mathcal Y}=\bigcup_{i=2}^n {\mathcal Y_i}$ with ${\mathcal Y}_i\cap{\mathcal Y}_j=\emptyset$ for $i\neq j$ and ${\mathcal Y}_m=\{y_{m1},\ldots,y_{mm}\}$, $2 \leq m \leq n$, and define $y_{m1} < \cdots < y_{mm}$ and $y_{mm} < y_{(m+1)1}$ for all $m \in \{2, \ldots, n-1\}$. 
For $m \in \{2, \ldots, n\}$, let $Y_m$ be the free $\Z$-module generated by ${\mathcal Y}_m$ and thus $L=L({\mathcal Y})=L(Y_2\oplus\ldots\oplus Y_n)$. For $r \in \{2, \ldots, n-1\}$, define $U_r,T_r$ and $\tilde{T}_r$ as in the previous subsection. Namely,  $U_{r} = Y_{r} \oplus \cdots \oplus Y_{n}$ with ${\mathcal U}_{r} = \bigcup_{\kappa = r}^{n}{\mathcal Y}_{\kappa}$ is a natural $\mathbb{Z}$-basis of $U_{r}$, $T_{r} = L(\Psi_{r}(W_{(r)}))$ and $\tilde{T}_{r} = \bigoplus_{s=2}^{r}T_{s}$. Furthermore, we define the map $\psi_{2,r}: [{\mathcal U}_{r+1},{\mathcal Y}_{r}] \rightarrow [U_{r+1},Y_{r}] \oplus L^{2}(U_{r+1})$ as follows:
\begin{flushleft} 
$(F1)~\psi_{2,r}([y_{mi},y_{ri}]) = [y_{mi},y_{ri}], ~~~~~~2 \leq r < m \leq n,$ \\
$(F2)~\psi_{2,r}([y_{mi},y_{rj}]) = [y_{mi},y_{rj}], \ \  2 \leq r < i \leq m \leq n, ~~1 \leq j \leq r,$ \\
$(F3)~\psi_{2,r}([y_{mi},y_{rj}]) = [y_{mi},y_{rj}]-[y_{mi},y_{mj}],\ \ 2 \leq r < m \leq n,~~ 1 \leq i \leq r < m \leq n,~~ 1 \leq j \leq r, j \neq i.$
\end{flushleft}

\begin{lemma}\label{le4}
For $n \geq 4$, $r \in \{3, \ldots, n-1\}$ and $j \in \{2, \ldots, r-1\}$, $$
[\psi_{2,r}([U_{r+1}, Y_{r}]), Y_{j}] \subseteq \tilde{T}_{r}.$$ 
\end{lemma}

\pf Let $n \geq 4$ and fix $r \in \{3, \ldots, n-1\}$ and $j \in \{2, \ldots, r-1\}$.
 Let 
$$
w = [\psi_{2,r}([y_{mk}, y_{rl}]),y_{jt}], ~{\rm with}~ r+1 \leq m \leq n.
$$
Note that $k \in \{1, \ldots, m\}$, $l \in \{1, \ldots, r\}$ and $t \in \{1,\ldots,j\}$. 
In what follows, we make extensive use of the definition of $\psi_{2,\mu}$ ($(F1), (F2)~ {\rm and}~ (F3)$) with $\mu \in \{2, \ldots, n-1\}$, the Jacobi identity of the form $[x,y,z] = [x,z,y] - [y,z,x]$ and the fact that $T_\nu$ is an ideal of $L(U_\nu)$ for all $\nu=2,\ldots,n-1$, a consequence of Proposition \ref{pr1}.  
We separate various cases. We give full details for the proof of the first case. Similar arguments may be applied to the remaining  cases. 
\begin{enumerate}
\item $k, l \leq j \leq r-1$.
\begin{enumerate}
\item $k \neq l$. Then, we separate the following cases. 
\begin{enumerate}

\item $t \neq k \leq j$ and $l \neq t \leq j$. Then, since $1 \leq k \neq l \leq r-1$, we have, by using $(F3)$,  
$$
\begin{array}{rll}
w & = & [\psi_{2,r}([y_{mk},y_{rl}]),y_{jt}], ~~~~r+1 \leq m \leq n, \\
& = & [y_{mk},y_{rl},y_{jt}] - [y_{mk},y_{ml},y_{jt}], ~~~~r+1 \leq m \leq n.
\end{array}
$$
By using the Jacobi identity in the form $[x,y,z] = [x,z,y] - [y,z,x]$, 
$$
\begin{array}{rll}
w & = & [y_{mk},y_{jt},y_{rl}] - [y_{rl},y_{jt},y_{mk}] \\
& & -[y_{mk},y_{jt},y_{ml}] + [y_{ml},y_{jt},y_{mk}].
\end{array}
$$
Since $1 \leq t \neq k \leq j$ and $1 \leq l \neq t \leq j$, we get, by using $(F3)$, 
$$
\begin{array}{rll}
w & = & [\psi_{2,j}([y_{mk},y_{jt}]),y_{rl}] + [y_{mk},y_{mt},y_{rl}] \\
& & -[\psi_{2,j}([y_{rl},y_{jt}]),y_{mk}] - [y_{rl},y_{rt},y_{mk}] \\
& & - [\psi_{2,j}([y_{mk},y_{jt}]),y_{ml}] - [y_{mk},y_{mt},y_{ml}] \\
& & +[\psi_{2,j}([y_{ml},y_{jt}]),y_{mk}] + [y_{ml},y_{mt},y_{rk}].
\end{array}
$$
Write
$$
\begin{array}{rll}
w_{1} & = & [\psi_{2,j}([y_{mk},y_{jt}]),y_{rl}] -[\psi_{2,j}([y_{rl},y_{jt}]),y_{mk}] \\
& & - [\psi_{2,j}([y_{mk},y_{jt}]),y_{ml}] +[\psi_{2,j}([y_{ml},y_{jt}]),y_{mk}].
\end{array}
$$
By using the Jacobi identity, we have 
$$
\begin{array}{rll}
w & = & w_{1} + [y_{mk},y_{rl},y_{mt}] - [y_{mt},y_{rl},y_{mk}] \\
& & + [y_{mk},y_{rl},y_{rt}] - [y_{mk},y_{rt},y_{rl}] \\
& & -[y_{mk},y_{mt},y_{ml}] + [y_{ml},y_{mt},y_{mk}].
\end{array}
$$
Since $1 \leq k \neq l < r$, $1 \leq t \neq l < r$ and $1 \leq k \neq t < r$, we obtain, by using $(F3)$,
$$
\begin{array}{rll}
w & = & w_{1} + [\psi_{2,r}([y_{mk},y_{rl}]),y_{mt}] + [y_{mk},y_{ml},y_{mt}] \\
& & -[\psi_{2,r}([y_{mt},y_{rl}]),y_{mk}] - [y_{mt},y_{ml},y_{mk}] \\
& & + [\psi_{2,r}([y_{mk},y_{rl}]),y_{rt}] + [y_{mk},y_{ml},y_{rt}] \\
& & - [\psi_{2,r}([y_{mk},y_{rt}]),y_{rl}] - [y_{mk},y_{mt},y_{rl}] \\
& & - [y_{mk},y_{mt},y_{ml}] + [y_{ml},y_{mt},y_{mk}].
\end{array}
$$
Write
$$
\begin{array}{rll}
w_{2} & = & [\psi_{2,r}([y_{mk},y_{rl}]),y_{mt}] -[\psi_{2,r}([y_{mt},y_{rl}]),y_{mk}] \\
& & + [\psi_{2,r}([y_{mk},y_{rl}]),y_{rt}] - [\psi_{2,r}([y_{mk},y_{rt}]),y_{rl}].
\end{array}
$$
We point out that 
$$
[y_{mk},y_{ml},y_{mt}] = [y_{mk},y_{mt},y_{ml}] - [y_{ml},y_{mt},y_{mk}]
$$
and so, 
$$
\begin{array}{rll}
w & = & w_{1} + w_{2} - [y_{mt},y_{ml},y_{mk}] + [y_{mk},y_{ml},y_{rt}] - [y_{mk},y_{mt},y_{rl}].
\end{array}
$$
By using the Jacobi identity,
$$
\begin{array}{rll}
w & = & w_{1} + w_{2} - [y_{mt},y_{ml},y_{mk}] + [y_{mk},y_{rt},y_{ml}] - [y_{ml},y_{rt},y_{mk}] \\
& & - [y_{mk},y_{rl},y_{mt}] + [y_{mt},y_{rl},y_{mk}].
\end{array}
$$
Since $1 \leq k \neq l < r$, $1 \leq t \neq l < r$ and $1 \leq k \neq t < r$, we obtain by using $(F3)$,
$$
\begin{array}{rll}
w & = & w_{1} + w_{2}-[y_{mt},y_{ml},y_{mk}]+[\psi_{2,r}([y_{mk},y_{rt}]),y_{ml}] + [y_{mk},y_{mt},y_{ml}] \\
& & -[\psi_{2,r}([y_{ml},y_{rt}]),y_{mk}] - [y_{ml},y_{mt},y_{mk}] -[\psi_{2,r}([y_{mk},y_{rl}]),y_{mt}] \\ 
& & - [y_{mk},y_{ml},y_{mt}] + [\psi_{2,r}([y_{mt},y_{rl}]),y_{mk}] + [y_{mt},y_{ml},y_{mk}].
\end{array}
$$
Finally, write 
$$
\begin{array}{rll}
w_{3} & = & [\psi_{2,r}([y_{mk},y_{rt}]),y_{ml}]-[\psi_{2,r}([y_{ml},y_{rt}]),y_{mk}] \\
& & -[\psi_{2,r}([y_{mk},y_{rl}]),y_{mt}] + [\psi_{2,r}([y_{mt},y_{rl}]),y_{mk}].
\end{array}
$$
Since $[y_{mk},y_{mt},y_{ml}] = [y_{mk},y_{ml},y_{mt}] - [y_{mt},y_{ml},y_{mk}]$, we obtain 
\begin{equation*}
\begin{split}
w = & w_{1} + w_{2} + w_{3} \\
= &[\psi_{2,j}([y_{mk},y_{jt}]),y_{rl}] - [\psi_{2,j}([y_{rl},y_{jt}]),y_{mk}] -  [\psi_{2,j}([y_{mk},y_{jt}]),y_{ml}] \\
&+ [\psi_{2,j}([y_{ml},y_{jt}]),y_{mk}] + 
[\psi_{2,r}([y_{mk},y_{rl}]),y_{rt}] - [\psi_{2,r}([y_{mk},y_{rt}]),y_{rl}]\\
&+ [\psi_{2,r}([y_{mk},y_{rt}]),y_{ml}] - [\psi_{2,r}([y_{ml},y_{rt}]),y_{mk}] \in  {\tilde T}_{r}.
\end{split}
\end{equation*}     

\item $1 \leq t \neq k \leq j$ and $1 \leq l = t \leq j$. Then, since $1 \leq k \neq t \leq j < r$ and $l = t$, we have, by using $(F3)$,
$$
\begin{array}{rll}
w & = & [\psi_{2,r}([y_{mk},y_{rt}]),y_{jt}], ~~~r+1 \leq m \leq n, \\
& = & [y_{mk},y_{rt},y_{jt}] - [y_{mk},y_{mt},y_{jt}], ~~r+1 \leq m \leq n.
\end{array}
$$
By using the Jacobi identity, 
$$
\begin{array}{rll}
w & = & [y_{mk},y_{jt},y_{rt}] - [y_{rt},y_{jt},y_{mk}] \\
& & -[y_{mk},y_{jt},y_{mt}] + [y_{mt},y_{jt},y_{mk}].
\end{array}
$$
Since $1 \leq t \neq k \leq j$, we get, by using $(F3)$ and $(F1)$, 
$$
\begin{array}{rll}
w & = & [\psi_{2,j}([y_{mk},y_{jt}]),y_{rt}] + [y_{mk},y_{mt},y_{rt}] \\
& & - [\psi_{2,j}([y_{rt},y_{jt}]),y_{mk}] - [\psi_{2,j}([y_{mk},y_{jt}]),y_{mt}] \\
& & -[y_{mk},y_{mt},y_{mt}] + [\psi_{2,j}([y_{mt},y_{jt}]),y_{mk}].
\end{array}
$$
Write 
$$
\begin{array}{rll}
w_{4} & = & [\psi_{2,j}([y_{mk},y_{jt}]),y_{rt}]- [\psi_{2,j}([y_{rt},y_{jt}]),y_{mk}] \\
& & - [\psi_{2,j}([y_{mk},y_{jt}]),y_{mt}]+ [\psi_{2,j}([y_{mt},y_{jt}]),y_{mk}]
\end{array}
$$
and by using the Jacobi identity, we have 
$$
w = w_{4} + [y_{mk},y_{rt},y_{mt}] - [y_{mt},y_{rt},y_{mk}] - [y_{mk},y_{mt},y_{mt}].
$$
Since $1 \leq t \neq k \leq j < r$, we obtain, by using $(F3)$,
$$
\begin{array}{rll}
w & = & w_{4} + [\psi_{2,r}([y_{mk},y_{rt}]),y_{mt}] - [\psi_{2,r}([y_{mt},y_{rt}]),y_{mk}].
\end{array}
$$
Thus,
\begin{equation*}
\begin{split}
w = &[\psi_{2,j}([y_{mk},y_{jt}]),y_{rt}] - [\psi_{2,j}([y_{rt},y_{jt}]),y_{mk}] -  [\psi_{2,j}([y_{mk},y_{jt}]),y_{mt}] \\
&+ [\psi_{2,j}([y_{mt},y_{jt}]),y_{mk}] + [\psi_{2,r}([y_{mk},y_{rt}]),y_{mt}] - [\psi_{2,r}([y_{mt},y_{rt}]),y_{mk}] \in  {\tilde T}_{r}.
\end{split}
\end{equation*}

\item $1 \leq t = k \leq j$ and $1 \leq l \neq t \leq j$. Then, since $1 \leq l \neq t \leq j < r$ and $k = t$, we have, by using $(F3)$,
$$
\begin{array}{rll}
w & = & [\psi_{2,r}([y_{mk},y_{rl}]),y_{jk}], ~~~r+1 \leq m \leq n, \\
& = & [y_{mk},y_{rl},y_{jk}] - [y_{mk},y_{ml},y_{jk}], ~~r+1 \leq m \leq n.
\end{array}
$$
By using the Jacobi identity, 
$$
\begin{array}{rll}
w & = & [y_{mk},y_{jk},y_{rl}] - [y_{rl},y_{jk},y_{mk}] \\
& & -[y_{mk},y_{jk},y_{ml}] + [y_{ml},y_{jk},y_{mk}].
\end{array}
$$
Since $1 \leq l \neq k \leq j$, we get, by using $(F1)$ and $(F3)$, 
$$
\begin{array}{rll}
w & = & [\psi_{2,j}([y_{mk},y_{jk}]),y_{rl}] - [\psi_{2,j}([y_{rl},y_{jk}]),y_{mk}] \\
& & - [y_{rl},y_{rk},y_{mk}] - [\psi_{2,j}([y_{mk},y_{jk}]),y_{ml}] \\
& & + [\psi_{2,j}([y_{ml},y_{jk}]),y_{mk}] + [y_{ml},y_{mk},y_{mk}].
\end{array}
$$
Write 
$$
\begin{array}{rll}
w_{5} & = & [\psi_{2,j}([y_{mk},y_{jk}]),y_{rl}] - [\psi_{2,j}([y_{rl},y_{jk}]),y_{mk}] \\
& & - [\psi_{2,j}([y_{mk},y_{jk}]),y_{ml}]+ [\psi_{2,j}([y_{ml},y_{jk}]),y_{mk}]. 
\end{array}
$$
Now, 
$$
w = w_{5} - [y_{rl},y_{rk},y_{mk}] + [y_{ml},y_{mk},y_{mk}].
$$
By using the Jacobi identity, 
$$
w = w_{5} - [y_{rl},y_{mk},y_{rk}] - [y_{mk},y_{rk},y_{rl}]  + [y_{ml},y_{mk},y_{mk}].
$$
Since $1 \leq l \neq k \leq j < r$, we obtain, by using $(F3)$,
$$
\begin{array}{rll}
w & = & w_{5} + [\psi_{2,r}([y_{mk},y_{rl}]),y_{rk}] + [y_{mk},y_{ml},y_{rk}] \\
& & - [\psi_{2,r}([y_{mk},y_{rk}]),y_{rl}] + [y_{ml},y_{mk},y_{mk}].
\end{array}
$$
Since $[y_{mk},y_{ml},y_{rk}] = [y_{mk},y_{rk},y_{ml}] - [y_{ml},y_{rk},y_{mk}]$ and $1 \leq l \neq k \leq j < r$, we obtain, by using $(F3)$, 
$$
\begin{array}{rll}
w & = & w_{5} + [\psi_{2,r}([y_{mk},y_{rl}]),y_{rk}] - [\psi_{2,r}([y_{mk},y_{rk}]),y_{rl}] \\
& & + [\psi_{2,r}([y_{mk},y_{rk}]),y_{ml}] - [\psi_{2,r}([y_{ml},y_{rk}]),y_{mk}].
\end{array}
$$ 
Thus,
\begin{equation*}
\begin{split}
w = &[\psi_{2,j}([y_{mk},y_{jk}]),y_{rl}] - [\psi_{2,j}([y_{rl},y_{jk}]),y_{mk}] -  [\psi_{2,j}([y_{mk},y_{jk}]),y_{ml}] \\
&+ [\psi_{2,j}([y_{ml},y_{jk}]),y_{mk}] + [\psi_{2,r}([y_{mk},y_{rl}]),y_{rk}] - [\psi_{2,r}([y_{mk},y_{rk}]),y_{rl}] \\
& + [\psi_{2,j}([y_{mk},y_{rk}]),y_{ml}] - [\psi_{2,j}([y_{ml},y_{rk}]),y_{mk}]  \in  {\tilde T}_{r}.
\end{split}
\end{equation*}   

\end{enumerate}

\item $k = l$. Then, 
$$
w = [\psi_{2,r}([y_{mk},y_{rk}]),y_{jt}] = [y_{mk},y_{rk},y_{jt}].
$$
By using the Jacobi identity,
$$
w = [y_{mk},y_{jt},y_{rk}] - [y_{rk},y_{jt},y_{mk}].
$$
Note that $1 \leq k, t \leq j$ and so, we separate two cases.
\begin{enumerate}

\item $k = t$. Then,
$$
\begin{array}{rll}
w & = & [y_{mk},y_{jk},y_{rk}] - [y_{rk},y_{jk},y_{mk}] \\
& = & [\psi_{2,j}([y_{mk},y_{jk}]),y_{rk}] - [\psi_{2,j}([y_{rk},y_{jk}]),y_{mk}] \in  {\tilde T}_{r}.
\end{array} 
$$  

\item $k \neq t \leq j$. Then, by using $(F3)$, 
$$
\begin{array}{rll}
w & = & [y_{mk},y_{jt}]),y_{rk}] - [y_{rk},y_{jt}, y_{mk}] \\
& = & [\psi_{2,j}([y_{mk},y_{jt}]),y_{rk}] + [y_{mk},y_{mt},y_{rk}] \\
& & - [\psi_{2,j}([y_{rk},y_{jt}]),y_{mk}] - [y_{rk},y_{rt},y_{mk}].
\end{array}
$$
Write
$$
w_{6} = [\psi_{2,j}([y_{mk},y_{jt}]),y_{rk}]- [\psi_{2,j}([y_{rk},y_{jt}]),y_{mk}].
$$
By using the Jacobi identity, $(F1)$ and $(F3)$, we get
$$
\begin{array}{rll}
w & = & w_{6} + [y_{mk},y_{rk},y_{mt}] - [y_{mt}, y_{rk},y_{mk}] \\
& & -[y_{rk},y_{mk},y_{rt}] - [y_{mk}, y_{rt},y_{rk}] \\
& = & w_{6} + [\psi_{2,r}([y_{mk},y_{rk}]),y_{mt}] - [\psi_{2,r}([y_{mt},y_{rk}]),y_{mk}] \\
& & -[y_{mt},y_{mk},y_{mk}] + [\psi_{2,r}([y_{mk},y_{rk}]),y_{rt}] \\
& & - [\psi_{2,r}([y_{mk},y_{rt}]),y_{rk}] - [y_{mk},y_{mt},y_{rk}].
\end{array}
$$ 
Write 
$$
\begin{array}{rll}
w_{7} & = & w_{6} + [\psi_{2,r}([y_{mk},y_{rk}]),y_{mt}] - [\psi_{2,r}([y_{mt},y_{rk}]),y_{mk}] \\
& & + [\psi_{2,r}([y_{mk},y_{rk}]),y_{rt}] - [\psi_{2,r}([y_{mk},y_{rt}]),y_{rk}]
\end{array}
$$
and so, 
$$
\begin{array}{rll}
w & = & w_{7} -[y_{mt},y_{mk},y_{mk}] - [y_{mk},y_{mt},y_{rk}] \\
({\rm Jacobi~identity}) & = & w_{7} - [y_{mt},y_{mk},y_{mk}] - [y_{mk},y_{rk},y_{mt}] + [y_{mt},y_{rk},y_{mk}] \\
& = & w_{7} - [\psi_{2,r}([y_{mk},y_{rk}]),y_{mt}]  + [\psi_{2,r}([y_{mt},y_{rk}]),y_{mk}]
\end{array}
$$
Thus,
\begin{equation*}
\begin{split}
w = &[\psi_{2,j}([y_{mk},y_{jt}]),y_{rk}] - [\psi_{2,j}([y_{rk},y_{jt}]),y_{mk}] + [\psi_{2,r}([y_{mk},y_{rk}]),y_{rt}] \\
&- [\psi_{2,r}([y_{mk},y_{rt}]),y_{rk}] \in {\tilde T}_{r}.
\end{split}
\end{equation*}

\end{enumerate} 
\end{enumerate}

\item $k \leq j \leq r-1$ and $j+1 \leq  l \leq r$. Then, 
\begin{equation*}
\begin{split}
w = &[\psi_{2,j}([y_{mk},y_{jt}]),y_{rl}] - [\psi_{2,j}([y_{rl},y_{jt}]),y_{mk}] -  [\psi_{2,j}([y_{mk},y_{jt}]),y_{ml}] \\
&+ [\psi_{2,j}([y_{ml},y_{jt}]),y_{mk}] - 
[\psi_{2,r}([y_{mt},y_{rl}]),y_{mk}] + [\psi_{2,r}([y_{mk},y_{rl}]),y_{mt}] \in  {\tilde T}_{r}.
\end{split}
\end{equation*}

\item $j+1 \leq k \leq r$ and $l \leq j \leq r$. Then, 
\begin{equation*}
\begin{split}
w = &[\psi_{2,j}([y_{mk},y_{jt}]),y_{rl}] - [\psi_{2,j}([y_{rl},y_{jt}]),y_{mk}] -  [\psi_{2,j}([y_{mk},y_{jt}]),y_{ml}] \\
&+ [\psi_{2,j}([y_{ml},y_{jt}]),y_{mk}] + 
[\psi_{2,r}([y_{mk},y_{rl}]),y_{rt}] - [\psi_{2,r}([y_{mk},y_{rt}]),y_{rl}] \\
&+  [\psi_{2,r}([y_{mk},y_{rt}]),y_{ml}] -  [\psi_{2,r}([y_{ml},y_{rt}]),y_{mk}] - [\psi_{2,r}([y_{mk},y_{rl}]),y_{mt}] \\
&+ [\psi_{2,r}([y_{mt},y_{rl}]),y_{mk}] \in  {\tilde T}_{r}.
\end{split}
\end{equation*}
\item $j+1 \leq k \leq r$ and $j+1 \leq l \leq r$. Then, 
\begin{equation*}
\begin{split}
w = &[\psi_{2,j}([y_{mk},y_{jt}]),y_{rl}] - [\psi_{2,j}([y_{rl},y_{jt}]),y_{mk}] -  [\psi_{2,j}([y_{mk},y_{jt}]),y_{ml}] \\
&+ [\psi_{2,j}([y_{ml},y_{jt}]),y_{mk}] \in  {\tilde T}_{r}.
\end{split}
\end{equation*}
\item $r+1 \leq k$
\begin{enumerate}
\item $j+1 \leq l$. Then, 
\begin{equation*}
w = [\psi_{2,j}([y_{mk},y_{jt}]),y_{rl}] - [\psi_{2,j}([y_{rl},y_{jt}]),y_{mk}] \in  {\tilde T}_{r}.
\end{equation*}
\item $l \leq j$. Then, 
\begin{equation*}
\begin{split}
w = &[\psi_{2,j}([y_{mk},y_{jt}]),y_{rl}] - [\psi_{2,j}([y_{rl},y_{jt}]),y_{mk}] + [\psi_{2,r}([y_{mk},y_{rl}]),y_{rt}] \\
&- [\psi_{2,r}([y_{mk},y_{rt}]),y_{rl}] \in  {\tilde T}_{r}.
\end{split}
\end{equation*} 
\end{enumerate}
\end{enumerate}

\begin{remark}\label{re1}\upshape{By the proof of Lemma \ref{le4}, we have, for $n \geq 4$, $r \in \{3, \ldots, n-1\}$ and $j \in \{2, \ldots, r-1\}$, 
$$
\begin{array}{rll}
[\psi_{2,r}([{\mathcal U}_{r+1},{\mathcal Y}_{r}]),{\mathcal Y}_{j}] & \subseteq & \left(\bigoplus_{a,b \geq 0 \atop a+b=1}[\psi_{2,j}([U_{j+1},Y_{j}]),~_{a}Y_{j},~_{b}U_{j+1}]\right) \oplus \\
& &  \left(\bigoplus_{a,b \geq 0 \atop a+b=1}[\psi_{2,r}([U_{r+1},Y_{r}]),~_{a}Y_{r},~_{b}U_{r+1}]\right).
\end{array}
$$
}
\end{remark}

\begin{lemma}\label{le5}
For $n \geq 4$, $r \in \{3, \ldots, n-1\}$, $j \in \{2, \ldots, r-1\}$ and $a, b \geq 0$, $
[\psi_{2,r}([U_{r+1}, Y_{r}]), ~_{a}Y_{r},~_{b}U_{r+1}, Y_{j}] \subseteq \tilde{T}_{r}$. 
\end{lemma}

\pf For $a = b = 0$, our claim follows from Lemma \ref{le4}. Thus, we assume that $a, b \geq 0$ with $a+b \geq 1$. Let  $
w_{a,b} = [u,y_{1}, \ldots, y_{a},y^{\prime}_{1}, \ldots, y^{\prime}_{b}]$  
with $u \in \psi_{2,r}([{\mathcal U}_{r+1},{\mathcal Y}_{r}])$, $y_{1}, \ldots, y_{a} \in {\mathcal Y}_{r}$, $y^{\prime}_{1}, \ldots, y^{\prime}_{b} \in {\mathcal{U}}_{r+1}$. 

Suppose that $a+b = 1$ and let $w = [u,y,z]$, where $u \in \psi_{2,r}([{\mathcal U}_{r+1},{\mathcal Y}_{r}])$, $y \in {\mathcal Y}_{s}$, $s\ge r$ and $z \in {\mathcal Y}_{j}$. By the Jacobi identity, $
w = [u,z,y] + [u,[y,z]]$. 
Since $j < r$, either $[y,z] = \psi_{2,j}([y,z])$ or $[y,z] = \psi_{2,j}([y,z]) + [y,y_{1}]$, with $y, y_{1} \in {\mathcal Y}_{s}$, $r\leq s$. Then, either 
$$
w = [u,z,y] + [u, \psi_{2,j}([y,z])]
$$
or
$$
w = [u,z,y] + [u,\psi_{2,j}([y,z])] + [u,y,y_{1}] - [u,y_{1},y].
$$
By Remark \ref{re1} and using the technique of Lemma \ref{le4} we get $w \in \tilde{T}_{r}$, since $\tilde{T}_{r}$ is a Lie subalgebra. 

So, we assume that $a+b \geq 2$. Then, by using repeatedly the Jacobi identity, we get 
$$
\begin{array}{rll}
[w_{a,b},z] & = & [[w_{a,b-1},z],y^{\prime}_{b}] + [w_{a,b-1},[y^{\prime}_{b},z]] \\
& = & [[w_{a,b-2},z],y^{\prime}_{b-1},y^{\prime}_{b}] + [w_{a,b-2},[y^{\prime}_{b-1},z],y^{\prime}_{b}] + [w_{a,b-1},[y^{\prime}_{b},z]] \\
& \vdots & \\
& = & [[u,z],y_{1}, \ldots, y_{a},y^{\prime}_{1}, \ldots, y^{\prime}_{b}] + [u,[y_{1},z],y_{2}, \ldots, y_{a},y^{\prime}_{1}, \ldots, y^{\prime}_{b}] +  \\
& & [w_{1,0},[y_{2},z],y_{3}, \ldots, y_{a},y^{\prime}_{1}, \ldots, y^{\prime}_{b}] + \ldots + [w_{a,b-2}, [y^{\prime}_{b-1},z],y^{\prime}_{b}] + \\ 
& & [w_{a,b-1},[y^{\prime}_{b},z]].
\end{array}
$$
As before, we obtain $[w_{a,b},z] \in \tilde{T}_{r}$. \qed    

\begin{lemma}\label{co3}
For $n \geq 4$, $r \in \{3, \ldots, n-1\}$, $j \in \{2, \ldots, r-1\}$, $a, b \geq 0$ and $s \geq 1$, $
[\psi_{2,r}([U_{r+1}, Y_{r}]), ~_{a}Y_{r},~_{b}U_{r+1}, ~_{s}Y_{j}] \subseteq \tilde{T}_{r}$. 
\end{lemma}

\pf For $a,b \geq 0$ and $s \geq 1$, we write 
$$
w(a,b,s) = [u, y_{1}, \ldots, y_{a},y^{\prime}_{1}, \ldots, y^{\prime}_{b}, z_{1}, \ldots, z_{s}],
$$
where $u \in \psi_{2,r}([{\mathcal U}_{r+1},{\mathcal Y}_{r}])$, $y_{1}, \ldots, y_{a} \in {\mathcal Y}_{r}$, $y^{\prime}_{1}, \ldots, y^{\prime}_{b} \in {\mathcal U}_{r+1}$ and $z_{1}, \ldots, z_{s} \in {\mathcal Y}_{j}$. It is enough to show that $w(a,b,s) \in \tilde{T}_{r}$. We use induction on $s$. By Lemma \ref{le5}, $w(a,b,1) \in \tilde{T}_{r}$. Suppose that $w(a,b,s-1) \in \tilde{T}_{r}$ and write 
$$
w(a,b,s-1) = u_{(2)} + u_{(3)} + \cdots + u_{(r)},
$$
where $u_{(j)} \in T_{j}$, $j = 2, \ldots, r$. Thus,
$$
w(a,b,s) = [u_{(2)},z] + [u_{(3)},z] + \cdots + [u_{(r)},z].
$$
Since $T_{2}$ is an ideal in $L$ and each $T_{j}$, $j \in \{2, \ldots,r\}$, is a free Lie algebra with a free generating set ${\mathcal C}_{(j)}$, we obtain from Lemma \ref{le5}, $w(a,b,s) \in \tilde{T}_{r}$ and so, we get the required result. \qed

\vskip .120 in

Let $J$ be the ideal of $L$ generated by the set ${\mathcal J}$ of all possible elements of the form
\begin{enumerate}
\item $[y_{mi},y_{ri}], ~~~~~~2 \leq r < m \leq n,$
\item $[y_{mi},y_{rj}], \ \  2 \leq r < i \leq m \leq n, ~~1 \leq j \leq r,$
\item $[y_{mi},y_{rj}]-[y_{mi},y_{mj}],\ \ 2 \leq r < m \leq n,~~ 1 \leq i \leq r < m \leq n,~~ 1 \leq j \leq r, j \neq i.$
\end{enumerate}

\begin{proposition}\label{pr3}
\begin{enumerate}
\item  For $n \geq 3$, $r \in \{2, \ldots, n-1\}$ and $j \in \{1, \ldots, r-1\}$, $
[T_{r}, L(Y_{j})] \subseteq \tilde{T}_{r}$.
 
\item For all $r \in \{2, \ldots, n-1\}$ with $n \geq 3$, $\tilde{T}_{r}$ is an ideal in $L$. Furthermore, $J = \tilde{T}_{n-1}$.
\end{enumerate}
\end{proposition}

\pf \begin{enumerate}
\item For $n = 3$, our claim is valid, since $T_{2}$ is an ideal in $L$. Thus, we may assume that $n \geq 4$. For $r \in \{3, \ldots, n-1\}$, $T_{r} = L({\mathcal C}_{(r)})$, where 
$$
{\mathcal C}_{(r),\kappa} = \bigcup_{a, b \geq 0 \atop a+b=\kappa -2}[\psi_{2,r}([{\mathcal U}_{r+1}, {\mathcal Y}_{r}]), ~_{a}{\mathcal Y}_{r}, ~_{b}{\mathcal U}_{r+1}]~ {\rm and}~ {\mathcal C}_{(r)} = \bigcup_{\kappa \geq 2}{\mathcal C}_{(r),\kappa}.
$$
Hence, for $r \in \{3, \ldots, n-1\}$ and $j \in \{1, \ldots, r-1\}$, it is enough to show that, for $a, b \geq 0$ and $s \geq 1$,
$$
[\psi_{2,r}([{\mathcal U}_{r+1},{\mathcal Y}_{r}]), ~_{a}{\mathcal Y}_{r}, ~_{b}{\mathcal U}_{r+1}, ~_{s}{\mathcal Y}_{j}] \subseteq \tilde{T}_{r},
$$
which is an immediate consequence of Lemma \ref{co3}. Therefore, 
$$[T_{r}, L(Y_{j})] \subseteq \tilde{T}_{r} \eqno(4)$$
for $r \in \{2, \ldots, n-1\}$ and $j \in \{1, \ldots, r-1\}$. 

\item Let $r \in \{2, \ldots, n-1\}$ with $n \geq 3$. To show that $\tilde{T}_{r}$ is an ideal in $L$, we induct on $r$. For $r = 2$, $\tilde{T}_{2} = T_{2}$, which is an ideal in $L$, and so, we suppose that our claim is valid for $k$ with $2 \leq k \leq n-2$ and $n \geq 4$. That is, $\tilde{T}_{k}$ is an ideal in $L$. By Proposition \ref{pr2} (2), 
$$
[T_{k+1},L] \subseteq \left(\sum_{j=2}^{k}[T_{k+1},L(Y_{j})]\right) + T_{k+1} + [T_{k+1},\tilde{T}_{k}].
$$
Since $\tilde{T}_{k+1} = T_{k+1} \oplus \tilde{T}_{k}$ and $\tilde{T}_{k}$ is an ideal in $L$, it is enough to show that $
[T_{k+1}, L(Y_{j})] \subseteq \tilde{T}_{k+1}$ 
for all $j \in \{2, \ldots, \kappa\}$. By the equation (4) (for $r = k+1$), we get $\tilde{T}_{r}$ is an ideal of $L$. Since ${\mathcal{J}} \subseteq \tilde{T}_{n-1}$, we obtain $J \subseteq \tilde{T}_{n-1}$. Since, for all $r \in \{2, \ldots, n-1\}$, $T_{r} = L({\mathcal C}_{(r)})$ and ${\mathcal C}_{(r)} \subseteq J$, we have $T_{r} \subseteq J$ for all $r$. Therefore, $\tilde{T}_{n-1} \subseteq J$ and so, $J = \tilde{T}_{n-1}$. \qed
\end{enumerate} 

\begin{theorem}\label{th1}
$L = (\bigoplus_{i=2}^{n}L(Y_{i})) \oplus J$.
\end{theorem}

\pf By Proposition \ref{pr2}, we have $
L = (\bigoplus_{i=2}^{n}L(Y_{i})) \oplus \tilde{T}_{n-1}$ 
and, by Proposition \ref{pr3}, we obtain the required result. \qed

\bigskip

Recall that for any $r \in \{2, \ldots, n-1\}$, $T_{r} = L(\Psi_{r}(W_{(r)}))$ and $
L(U_{r}) = L(Y_{r}) \oplus L(U_{r+1}) \oplus T_{r}$. 
By Corollary \ref{co2}, for any $m \geq 2$, 
$$
L^{m}(U_{r}) = L^{m}(Y_{r}) \oplus L^{m}(U_{r+1}) \oplus L^{m}_{\rm grad}(\Psi_{r}(W_{(r)})). \eqno(5)
$$
By Proposition \ref{pr2} (1) and the equation (5), we have the following result.

\begin{corollary}
With the above notation, for any $m \geq 2$, $
L^{m} = (\bigoplus_{i=2}^{n}L^{m}(Y_{i})) \oplus (\bigoplus_{j=2}^{n-1}L^{m}_{\rm grad}(\Psi_{j}(W_{(j)})))$. 
In particular, for any $m \geq 2$, $
J^{m} = \bigoplus_{j=2}^{n-1}L^{m}_{\rm grad}(\Psi_{j}(W_{(j)}))$.  
\end{corollary}
    
\subsection{A presentation of ${\rm gr}(I_{n})$}

Let $F=F({\mathcal Y})$ be the free group on the finite set ${\mathcal Y}$ such that ${\mathcal Y}$ is decomposed into a disjoint union ${\mathcal Y}={\mathcal Y}_2\cup{\mathcal Y}_3\cup\ldots\cup{\mathcal Y}_n$ with $n\ge 3$, of subsets of the form ${\mathcal Y}_m=\{y_{m1},\ldots,y_{mm}\}$, with $2\le m\le n$, and such that
$$y_{21}<y_{22}<y_{31}<y_{32}<y_{33}<\ldots<y_{n1}<\ldots<y_{nn}.$$ 
It is well known that $\gr(F)$ is a free Lie algebra of finite rank $|{\mathcal Y}|$ with a free generating set $\{y_{ij}F' \mid y_{ij}\in {\mathcal Y}\}$. Both $L$ and $\gr(F)$ are free Lie algebras of the same rank and they are isomorphic to each other as Lie algebras in a natural way, 
and from now on, we identify $L$ and $\gr(F)$. For all $c\ge 1$, we have $L^c=\g_c(F)/\g_{c+1}(F)$.  

Let ${\mathcal R}$ be the subset of $\gamma_{2}(F) = F^{\prime}$ consisting of all possible elements of the form
\begin{enumerate}
\item $[y_{mi},y_{ri}], ~~~~~~2 \leq r < m \leq n,$
\item $[y_{mi},y_{rj}], \ \  2 \leq r < i \leq m \leq n, ~~1 \leq j \leq r,$
\item $[y_{mi},y_{rj}][y_{mi},y_{mj}]^{-1},\ \ 2 \leq r < m \leq n,~~ 1 \leq i \leq r < m \leq n,~~ 1 \leq j \leq r, j \neq i.$
\end{enumerate}
Let $N = {\mathcal R}^{F}$ be the normal closure of ${\mathcal R}$ in $F$. Thus, $N$ is generated by the set $\{r^{g}= g^{-1}rg: r \in {\mathcal R}, g \in F\}$. By using the presentation of $I_{n}$, the group $I_{n}$ may be identified with the quotient group $F/N$. Since $r \in \gamma_{2}(F) \setminus \gamma_{3}(F)$ for all $r \in {\mathcal R}$, we have $N \subseteq F^{\prime}$ and so, $
I_{n}/I^{\prime}_{n} \cong F/F^{\prime}$. Let $
\widetilde{{\mathcal R}} = \{[r,g]: r \in {\mathcal R}, g \in F \setminus \{1\}\} \subseteq \gamma_{3}(F)$. Since $N\gamma_{3}(F)/\gamma_{3}(F)$ is generated by the set $\{r\gamma_{3}(F): r \in {\mathcal R}\}$, we have $
N\gamma_{3}(F)/\gamma_{3}(F) = J^{2}$,  
where $J^{2} = \bigoplus_{j=2}^{n-1}\Psi_{j}(W_{(j),2})$. Since $L^{2} = (\bigoplus_{i=2}^{n}L^{2}(Y_{i})) \oplus J^{2}$, we get the set $\{r\gamma_{3}(F): r \in {\mathcal R}\}$ is $\mathbb{Z}$-linear independent and so, $
{\mathcal R} \cap \widetilde{\mathcal R} = \emptyset$. 
Therefore, $N$ is generated by the disjoint union ${\mathcal R} \cup \widetilde{\mathcal R}$.

For a positive integer $d$, let $N_{d} = N \cap \gamma_{d}(F)$.
Note that for $d \leq 2$, we have $N_{d} = N$. Also, for $d \geq
2$, 
$$
\begin{array}{lll}
N_{d+1} &  = & N \cap \gamma_{d+1}(F) \\
& = & (N \cap \gamma_{d}(F)) \cap \gamma_{d+1}(F) \\ 
& = & N_{d} \cap \gamma_{d+1}(F).
\end{array}
$$
Since $[g_{1},
\ldots, g_{\kappa}]^{f} = [g_{1}^{f}, \ldots, g_{\kappa}^{f}]$ for
all $f \in F$, and $N$ is normal, we obtain $\{N_{d}\}_{d \geq 2}$
is a normal (descending) series of $N$. Clearly, each $N_{d}$ is
normal in $F$. Since $[N_{\kappa}, N_{\ell}] \subseteq N_{\kappa +
\ell}$ for all $\kappa, \ell \geq 2$, we have $\{N_{d}\}_{d \geq
2}$ is a central series of $N$. Define $
{\mathcal I}_{d}(N) = N_{d}
\gamma_{d+1}(F)/\gamma_{d+1}(F) \leq \gamma_{d}(F)/\gamma_{d+1}(F)$.  
Since $N \subseteq F^{\prime}$, we get ${\mathcal
I}_{1}(N) = \{F^{\prime}\}$. 
Since $N$ is a normal subgroup of $F$, we have the Lie subalgebra $
{\mathcal I}(N) = \bigoplus_{d \geq 2}{\mathcal I}_{d}(N)$ of ${\rm gr}(F)$  is an
ideal of ${\rm gr}(F)$ (see \cite{laz}). It is easily verified that
${\mathcal I}_{d}(N) \cong N_{d}/N_{d+1}$ as $\mathbb{Z}$-modules. By
our definitions, identifications and the above discussion, ${\mathcal
I}_{2}(N) = N\gamma_{3}(F)/\gamma_{3}(F) = J^{2}$. Since ${\mathcal
I}_{2}(N) = J^{2}$, ${\mathcal I}(N)$ is an ideal of ${\rm gr}(F) = L$ and $J$ is generated 
by ${\mathcal J}$, we obtain $J \subseteq {\mathcal I}(N)$. Since $J$ is a homogeneous ideal, we get, for all $d \geq 2$, $J^{d} \subseteq {\mathcal I}_{d}(N)$. For a positive integer $m$, with $m \geq 2$, we write $
(L^{m})^{*} = \bigoplus_{i=2}^{n}L^{m}(Y_{i})$. 
Denote $
{\mathcal V} = {\mathcal J} = \bigcup_{r=2}^{n-1}\psi_{2,r}([{\mathcal U}_{r+1}, {\mathcal Y}_{r}])$ 
and $
{\mathcal V}^{*} = \bigcup_{m=2}^{n}\{[y_{m\kappa}, y_{m \ell}]: 1 \leq \ell < \kappa \leq m\}$.  
We point out that ${\mathcal V}^{*}$ is a natural $\mathbb{Z}$-basis of $(L^{2})^{*}$ and ${\mathcal R}\gamma_{3}(F)/\gamma_{3}(F) = {\mathcal V} = {\mathcal J}$, 
which is a natural $\mathbb{Z}$-basis
of $J^{2}$.

\begin{proposition}\label{pr4}
For all positive integers $c$, ${\mathcal
I}_{c+2}(N) = J^{c+2}$.
\end{proposition}

\pf Let $c = 1$. We denote by $N_{2,1}$ the subgroup of $N_{2} = N$ generated by the set $\mathcal R$  and let $N_{2,2}$ be the normal closure of $\widetilde{\mathcal R}$ in $F$. Since $N_{2} = N = N_{2,1}N_{2,2}$ and $N_{2,2} \subseteq \gamma_{3}(F)$, by using the modular law, we get 
$$
\begin{array}{rll}
N_{3} & = & N_{2} \cap \gamma_{3}(F) \\
& = & (N_{2,1}N_{2,2}) \cap \gamma_{3}(F) \\
& = & (N_{2,1} \cap \gamma_{3}(F)) N_{2,2}.
\end{array}
$$
An element $w$ of $N_{2,1}$ is written as a product of elements in ${\mathcal R}$ (in an abelian form) and an element $u \in N^{\prime}_{2,1}$. Since ${\mathcal I}_{2}(N) = J^{2}$ and ${\mathcal R}\gamma_{3}(F)/\gamma_{3}(F)$ is a $\mathbb{Z}$-basis of ${\mathcal I}_{2}(N)$, we have  $N_{2,1} \cap \gamma_{3}(F)$ is generated by elements which belong to $N_{2,2} \cap \gamma_{4}(F)$. That is, $N_{2,1} \cap \gamma_{3}(F) \subseteq N_{2,2}$. Therefore, $N_{3} = N_{2,2}$ is generated by the set $\{[r,g]^{f}: r \in {\mathcal R}, g \in F \setminus \{1\}, f \in F\}$ and so, ${\mathcal I}_{3}(N) = J^{3}$. 

Suppose that $c \geq 2$. Now, $I_n=H_n\rtimes H_{n-1}\rtimes \ldots \rtimes H_2$ and so, by induction on $n$ and a result of Falk and Randell \cite[Theorem (3.1)]{fr}, we have, for all $d \geq 3$, $
\gamma_{d}(I_{n})/\gamma_{d+1}(I_{n}) \cong \bigoplus_{i=2}^{n} \gamma_{d}(H_{i})/\gamma_{d+1}(H_{i})$. 
Since $H_{i}$ is a free group of rank $i$, $i \in \{2, \ldots, n\}$, we have $
\gamma_{d}(H_{i})/\gamma_{d+1}(H_{i}) \cong L^{d}(Y_{i})$. 
By Proposition \ref{pr3} (2), Lemma \ref{l3} and Theorem \ref{th1},  
$$
\gamma_{d}(I_{n})/\gamma_{d+1}(I_{n}) \cong L^{d}/J^{d} \cong (L^{d})^{*} = \bigoplus_{i=2}^{n}L^{d}(Y_{i}).
$$
Hence, $
{\rm rank}(\gamma_{d}(I_{n})/\gamma_{d+1}(I_{n})) = {\rm rank}(L^{d})^{*}$. 
Now, 
$$
\gamma_{d}(I_{n})/\gamma_{d+1}(I_{n}) \cong \gamma_{d}(F)/(\gamma_{d}(F) \cap (\gamma_{d+1}(F) N)).
$$ 
Since $\gamma_{d+1}(F) \subseteq \gamma_{d}(F)$, we have, by the modular law, 
$$
\gamma_{d}(F)/(\gamma_{d}(F) \cap (\gamma_{d+1}(F) N)) = \gamma_{d}(F)/(\gamma_{d+1}(F) N_{d}).
$$
But, $
\gamma_{d}(F)/(\gamma_{d+1}(F) N_{d}) \cong (\gamma_{d}(F)/\gamma_{d+1}(F))/{\mathcal I}_{d}(N)$. 
Therefore, for all $d \geq 3$, $
\gamma_{d}(I_{n})/\gamma_{d+1}(I_{n}) \cong (\gamma_{d}(F)/\gamma_{d+1}(F))/{\mathcal I}_{d}(N)$ 
and so, 
$$
(\gamma_{d}(F)/\gamma_{d+1}(F))/{\mathcal I}_{d}(N) \cong L^{d}/J^{d} \cong (L^{d})^{*}. \eqno(6) 
$$
Since $J^{c} \subseteq {\mathcal I}_{c}(N)$ for all $c \geq 2$, $(L^{d})^{*}$ is torsion-free for all $d \geq 3$ and by our definitions and identifications, we conclude from $(6)$ that ${\mathcal I}_{d}(N) = J^{d}$. Thus, we obtain the required result. \qed 

\begin{corollary}
${\mathcal I}(N) = J$.
\end{corollary}

\pf Since $J = \bigoplus_{d \geq 2}J^{d}$ and ${\mathcal I}_{2}(N) =
J^{2}$, we have from Proposition \ref{pr4} that ${\mathcal I}(N) = J$.
\qed

\vskip .120 in

So, we are now able to prove the following.

\begin{theorem}\label{3} $L/J \cong {\rm gr}(I_{n})$ as Lie
algebras. 
\end{theorem}

\pf Recall that ${\rm gr}(I_n)=\bigoplus_{c\ge 1}\gamma_c(I_n)/\gamma_c(I_{n+1})$. Since $I_n/I_n'\cong F/NF'=F/F'\cong \Z^k$ with $k=(n-1)(n+2)/2$, we have ${\rm gr}(I_n)$ is generated as a Lie algebra by the set $\{\alpha_{ij} : 2\le i \le n, 1\le j\le i\}$ with $\alpha_{ij}=y_{ij}I^{\prime}_n$. Since $L$ is a free Lie algebra of rank $k$ with a free generating set $\{y_{ij} : 2\le i\le n, 1\le j\le i\}$, the map $\xi:L\ar {\rm gr}(I_n)$ with $\xi(y_{ij})=\alpha_{ij}$ extends uniquely to a Lie algebra epimorphism. Hence $L/\ker\xi\cong{\rm gr}(I_n)$ as Lie algebras. By definition, $J\subseteq\ker\xi$ and so, $\xi$ induces a Lie algebra epimorphism $\overline{\xi}:L/J\ar {\rm gr}(I_n)$. In particular, $\overline{\xi}(y_{ij}+J)=\alpha_{ij}$ for all $i,j$. Moreover, $\overline{\xi}$ induces a $\Z$-linear mapping $\overline{\xi}_c:(L^c+J)/J\ar \gamma_c(I_n)/\gamma_{c+1}(I_n)$. For $c\ge 2$, and in view of the proof of Proposition \ref{pr4}, we have ${\rm rank}(\g_c(I_n)/\g_{c+1}(I_n))={\rm rank}(L^c)^*$. Since $J = \bigoplus_{m \geq 2}J^{m}$, we have
$(L^{c} + J)/J \cong L^{c}/(L^{c} \cap J) = L^{c}/J^{c} \cong
(L^{c})^{*}$ (by Theorem \ref{th1} and Lemma \ref{l3}) and so, we obtain ${\rm
ker}\overline{\xi}_{c}$ is torsion-free. Since ${\rm
rank}(\gamma_{c}(I_n)/\gamma_{c+1}(I_n)) = {\rm
rank}(L^{c})^{*}$, we have ${\rm ker}\overline{\xi}_{c} = \{1\}$
and so, $\overline{\xi}_{c}$ is an isomorphism. Since
$\overline{\xi}$ is an epimorphism and each $\overline{\xi}_{c}$ is
isomorphism, we have $\overline{\xi}$ is an isomorphism. Hence, $L/J
\cong {\rm gr}(I_n)$ as Lie algebras. \qed

\begin{corollary}
A presentation of $\gr(I_n)$ is given by generators $y_{ij}$, with $i \in\{2,\ldots,n\}$, $j\in\{1,\ldots,n\}$ and $i\ge j$, together with the following relations:
\begin{enumerate}
\item $[y_{mi},y_{ri}], ~~~~~~2 \leq r < m \leq n,$
\item $[y_{mi},y_{rj}], \ \  2 \leq r < i \leq m \leq n, ~~1 \leq j \leq r,$
\item $[y_{mi},y_{rj}]-[y_{mi},y_{mj}],\ \ 2 \leq r < m \leq n,~~ 1 \leq i \leq r < m \leq n,~~ 1 \leq j \leq r, j \neq i.$
\end{enumerate} 
\end{corollary}

\section{A natural embedding}
For positive integers $n$ and $c$, with $n, c \geq 2$, let $F_n$ be a free group of rank $n$ and $F_{n,c-1}=F_n/\gamma_c(F_n)$ be the free nilpotent group of rank $n$ and class $c-1$. The natural epimorphism from $F_n$ onto $F_{n,c-1}$ induces a group homomorphism $\pi_{n,c-1}:\aut(F_n)\ar\aut(F_{n,c-1})$ 
with kernel ${\rm I}_c{\rm A}(F_n)$. For $c=2$ we write ${\rm I}_2{\rm A}(F_n)={\rm IA}(F_n)$ for the $IA$-automorphisms of $F_n$. It is well known that, for $t, s \geq
2$, $[{\rm I}_{t}{\rm A}(F_{n}), {\rm I}_{s}{\rm A}(F_{n})] \subseteq {\rm
I}_{t+s-1}{\rm A}(F_{n})$.  
Since $F_{n}$ is residually nilpotent, we have $\bigcap_{c \geq 2}{\rm I}_{c}{\rm A}(F_{n}) = \{1\}$.  
For $r \geq 2$, we write ${\mathcal L}^{r}({\rm I}{\rm A}(F_{n}))
= {\rm I}_{r}{\rm A}(F_{n})/{\rm I}_{r+1}{\rm
A}(F_{n})$. It is known that ${\mathcal L}^{r}({\rm I}{\rm A}(F_{n}))$ is a free abelian group for all 
$r \geq 2$ (see, for example, \cite[Section 1]{mps}). Form the (restricted) direct sum of the free 
abelian groups ${\mathcal L}^{r}({\rm I}{\rm A}(F_{n}))$ and
denote it by $
{\mathcal L}({\rm IA}(F_{n})) = \bigoplus_{r \geq 2}{\mathcal L}^{r}({\rm I}{\rm A}(F_{n}))$.  
It has the structure of a graded Lie algebra with ${\mathcal L}^{r}({\rm IA}(F_{n}))$ as component of degree $r-1$ in the grading and Lie multiplication given by $
[\phi {\rm I}_{j+1}{\rm A}(F_{n}), \psi {\rm I}_{\kappa+1}{\rm A}(F_{n})] = (\phi^{-1}\psi^{-1}\phi\psi){\rm I}_{j+\kappa}{\rm A}(F_{n})$,  
for all $\phi \in {\rm I}_{j}{\rm A}(F_{n})$, $\psi \in {\rm I}_{\kappa}{\rm A}(F_{n})$ with $j, \kappa \geq 2$. Multiplication is then extended to
${\mathcal L}({\rm IA}(F_{n}))$ by linearity. The above Lie algebra is usually called \emph{the
Andreadakis-Johnson  Lie algebra of} ${\rm IA}(F_{n})$. We point out that, for any positive integer
$c$, $ \gamma_{c}({\rm IA}(F_{n})) \subseteq {\rm I}_{c+1}{\rm
A}(F_{n})$. 

Let $G$ be a finitely generated subgroup of ${\rm IA}(F_{n})$ with $G/G^{\prime}$ torsion-free.  For
a positive integer $q$, let ${\mathcal L}^{q}_{1}(G) =
\gamma_{q}(G)({\rm I}_{q+2}{\rm A}(F_{n}))/{\rm I}_{q+2}{\rm
A}(F_{n})$. Form the
(restricted) direct sum of free abelian groups $
{\mathcal L}_{1}(G) = \bigoplus_{q \geq 1} {\mathcal L}^{q}_{1}(G)$.  
It is easily verified that ${\mathcal
L}_{1}(G)$ is a Lie subalgebra of ${\mathcal
L}({\rm IA}(F_{n}))$. Furthermore, if $\{y_{1}G^{\prime}, \ldots, y_{m}G^{\prime}\}$ is a $\mathbb{Z}$-basis for $G/G^{\prime}$, then ${\mathcal L}_{1}(G)$ is generated as Lie algebra by the set $$\{y_{1} ({\rm I}_{3}{\rm A}(F_{n})), \ldots, y_{m}({\rm I}_{3}{\rm A}(F_{n}))\}.$$ By \emph{a natural embedding} of ${\rm gr}(G)$ into ${\mathcal L}({\rm IA}(F_{n}))$, we mean that there exists a Lie algebra isomorphism $\zeta$ from ${\rm gr}(G)$ onto ${\mathcal L}_{1}(G)$ satisfying the conditions $\zeta(y_{i}G^{\prime}) = y_{i}({\rm I}_{3}{\rm A}(F_{n}))$, $i = 1, \ldots, m$. In this case, we also say that ${\rm gr}(G)$ is naturally isomorphic to ${\mathcal L}_{1}(G)$.

In this section we prove that there is a natural embedding of ${\rm gr}(I_n)$ into ${\mathcal L}({\rm IA}(F_n))$.

\begin{theorem}\label{kostas}
For a positive integer $n\ge 2$, ${\rm gr}(I_n)$ is naturally isomorphic to ${\mathcal L}_{1}(I_n)$ as Lie algebra.
\end{theorem}

\pf Let $n = 2$. Since $I_{2} = H_{2} = {\rm Inn}(F_{2}) \cong F_{2}$, our claim follows from \cite[Proof of Proposition 3]{mps}. Thus, we assume that $n \geq 3$. Recall that $I_n=(\ldots(H_n\rtimes H_{n-1})\rtimes\ldots)\rtimes H_2$, where each $H_i$ is a free group of finite rank $i$ with a free generating set $\{y_{i1}, \ldots, y_{ii}\}$. Thus, ${\rm gr}(H_{i}) \cong {\rm gr}(F_{i})$ as Lie algebras. By \cite[Proof of Proposition 3]{mps}, we have ${\rm gr}(H_{i}) \cong {\mathcal L}_{1}(H_{i})$ by a Lie algebra isomorphism $\zeta_{i}$ satisfying the conditions $\zeta_{i}(y_{ik}H^{\prime}_{i}) = y_{ik}({\rm I}_{3}{\rm A}(F_{i}))$, with $k = 1, \ldots, i$, for all $i = 2, \ldots, n$. Hence, ${\mathcal L}_{1}(H_{i})$ is a free Lie algebra of rank $i$. Since $y_{i1} \notin {\rm I}_{3}{\rm A}(F_{i})$, we have ${\mathcal L}_{1}(H_{i})$ is a non-trivial subalgebra of ${\mathcal L}({\rm IA}(F_{i}))$. In fact, ${\mathcal L}_{1}(H_{i})$ is generated by the set $\{y_{ik}({\rm I}_{3}{\rm A}(F_{i})) : k = 1, \ldots, i\}$. For $h \in F_{i}$, we write $\tau_{h}$ for the inner automorphism of $F_{i}$ defined by $\tau_{h}(x) = hxh^{-1}$ for all $x \in F_{i}$. Since $[\tau_{g},\phi] = \tau_{g^{-1}\phi^{-1}(g)}$ for all $\phi \in {\rm Aut}(F_{i})$ and $g \in F_{i}$, we have ${\mathcal L}_{1}(H_{i})$ is an ideal in ${\mathcal L}({\rm IA}(F_{i}))$ for all $i \in \{2, \ldots, n\}$. 

Next, we claim that ${\mathcal L}_{1}(I_n)=\bigoplus_{i=2}^{n}{\mathcal L}_{1}(H_i)$.  Fix $m \in \{3, \ldots, n\}$. First, we observe that $I_{m} = H_{m} \rtimes I_{m-1}$ and claim that ${\mathcal L}_{1}(I_{m})$ is additively equal to the direct sum of ${\mathcal L}_{1}(I_{m-1})$ and ${\mathcal L}_{1}(H_{m})$. Since ${\mathcal L}_{1}(H_{m})$ is an ideal in ${\mathcal L}({\rm IA}(F_{m}))$, we have ${\mathcal L}_{1}(I_{m-1}) + {\mathcal L}_{1}(H_{m})$ is a Lie subalgebra of ${\mathcal L}_{1}(I_{m})$. Let $w \in {\mathcal L}_{1}(I_{m-1}) \cap {\mathcal L}_{1}(H_{m})$. In the next few lines, we adopt the arguments given in the proof of \cite[Proposition 4]{mps}. Since both ${\mathcal L}_{1}(I_{m-1})$ and ${\mathcal L}_{1}(H_{m})$ are graded Lie algebras, we may assume that $w \in {\mathcal L}^{d}_{1}(I_{m-1})\cap {\mathcal L}^{d}_{1}(H_{m})$ for some $d$. Thus, there are $u \in \gamma_{d}(I_{m-1})$ and $v \in \gamma_{d}(H_{m})$ such that $w = u{\rm I}_{d+2}{\rm A}(F_{m}) = v{\rm I}_{d+2}{\rm A}(F_{m})$. To get a contradiction, we assume that $u, v \notin {\rm I}_{d+2}{\rm A}(F_{m})$. Therefore, $v \in \gamma_{d}(H_{m}) \setminus \gamma_{d+1}(H_{m})$ and so, there exists $g \in \gamma_{d}(F_{m}) \setminus \gamma_{d+1}(F_{m})$ such that $v = \tau_{g}\rho$, where $\rho \in \gamma_{d+1}(H_{m})$. Then, we have $u{\rm I}_{d+2}{\rm A}(F_{m}) =\tau_g\rho {\rm I}_{d+2}(F_m)=\tau_g {\rm I}_{d+2}(F_m)$, since $\rho \in \gamma_{d+1}(H_{m}) \subseteq {\rm I}_{d+2}{\rm A}(F_{m})$. So,  $u^{-1}\tau_{g} \in {\rm I}_{d+2}{\rm A}(F_{m})$. Hence, $u^{-1} \tau_{g}(x_{m}) = x_{m} h$, where $h \in \gamma_{d+2}(F_{m})$, by the definition of ${\rm I}_{d+2}{\rm A}(F_{m})$.
Since $u^{-1}$ fixes $x_{m}$, we have 
$$
\begin{array}{rll}
h = x^{-1}_{m}(u^{-1}\tau_{g}(x_{m})) & = & x^{-1}_{m}u^{-1}(g) x_{m} u^{-1}(g^{-1}) \\
& = & x^{-1}_{m}\tau_{u^{-1}(g)}(x_{m}) \\
& = & [x_{m}, u^{-1}(g^{-1})] \in \gamma_{d+2}(F_{m}).
\end{array}
$$
Since $I_{m-1} \subseteq {\rm IA}(F_{m})$, we obtain $\gamma_{d}(I_{m-1}) \subseteq \gamma_{d}({\rm IA}(F_{m})) \subseteq {\rm I}_{d+1}{\rm A}(F_{m})$ and so, $u^{-1}(x_{j}) = x_{j}z_{j}$, where $z_{j} \in \gamma_{d+1}(F_{m})$. Since $\gamma_{d}(F_{m})$ is a fully invariant subgroup of $F_{m}$, we get $u^{-1}(g^{-1}) = g^{-1}g_{1}$
 with $g_{1} \in \gamma_{d+1}(F_{m})$. Moreover, since 
$$
\begin{array}{rll}
[x_{m},u^{-1}(g^{-1}) & = & [x_{m},g^{-1}g_{1}] \\
& = & [x_{m},g_{1}][x_{m},g^{-1}][x_{m},g^{-1},g_{1}]
\end{array}
$$
and $[x_{m},g^{-1}] = ([x_{m},g]^{-1})^{g^{-1}}$, we have $[g,x_{m}] \in \gamma_{d+2}(F_{m})$. 
Since ${\rm gr}(F_{m})$ is a free Lie algebra of rank $m$ with a free generating set $\{x_{i}F^{\prime}_{m}: i = 1, \ldots, m\}$ and $\gamma_{d}(F_{m})/\gamma_{d+1}(F_{m})$ is the $d$th homogeneous component of ${\rm gr}(F_{m})$, we obtain $[g,x_{m}] \in \gamma_{d+1}(F_{m}) \setminus \gamma_{d+2}(F_{m})$ (see, for example, 
\cite[Theorem 5.10]{mks}) which is a contradiction. Therefore, ${\mathcal L}_{1}(I_{m-1}) 
\cap {\mathcal L}_{1}(H_{m}) = \{0\}$ for all $m \in \{3, \ldots, n\}$. Since $I_{m} = H_{m} \rtimes I_{m-1}$ and the action of $I_{m-1}$ on $H_{m}$ is by conjugation, we have, by a result of Falk and Randell \cite[Proof of Theorem (3.1)]{fr}, $\gamma_{q}(I_{m}) = \gamma_{q}(H_{m}) \rtimes \gamma_{q}(I_{m-1})$ for all $q$. Hence, ${\mathcal L}^{q}_{1}(I_{m}) = {\mathcal L}^{q}_{1}(I_{m-1}) + {\mathcal L}^{q}_{1}(H_{m})$ and so, ${\mathcal L}^{q}_{1}(I_{m}) = {\mathcal L}^{q}_{1}(I_{m-1}) \oplus {\mathcal L}^{q}_{1}(H_{m})$ for all $q$. Therefore, ${\mathcal L}_{1}(I_{m}) = {\mathcal L}_{1}(I_{m-1}) \oplus {\mathcal L}_{1}(H_{m})$ for all $m \in \{3, \ldots, n\}$ and so, 
$$
\begin{array}{rll}
{\mathcal L}_{1}(I_{n}) & = & {\mathcal L}_{1}(I_{n-1}) \oplus {\mathcal L}_{1}(H_{n}) \\
& = & {\mathcal L}_{1}(I_{n-2}) \oplus {\mathcal L}_{1}(H_{n-1}) \oplus {\mathcal L}_{1}(H_{n}) \\
& \vdots & \\
& = & {\mathcal L}_{1}(H_{2}) \oplus {\mathcal L}_{1}(H_{3}) \oplus \cdots \oplus {\mathcal L}_{1}(H_{n}).
\end{array}
$$
Recall that, for $m \in \{2, \ldots, n\}$, we have ${\rm gr}(H_{m})\cong{\mathcal L}_{1}(H_{m})$ by an isomorphism $\zeta_{m}$ satisfying the conditions $\zeta_{m}(y_{mk}H^{\prime}_{m}) = y_{mk}{\rm I}_{3}{\rm A}(F_{m})$ for all $k \in \{1, \ldots, m\}$. By a result of Andreadakis \cite[Section 6]{a}, $\gamma_{c}(H_{m}) = H_{m}\cap {\rm I}_{c+1}{\rm A}(F_{m})$ for all $c$ and $m \in \{2, \ldots, n\}$. Since 
$$
\begin{array}{rll}
{\mathcal L}_{1}^{c}(H_{m}) & = & \gamma_{c}(H_{m}){\rm I}_{c+2}{\rm A}(F_{m})/{\rm I}_{c+2}{\rm A}(F_{m}) \\
& \cong & \gamma_{c}(H_{m})/(\gamma_{c}(H_{m}) \cap {\rm I}_{c+2}{\rm A}(F_{m})) \\
& = & \gamma_{c}(H_{m})/\gamma_{c+1}(H_{m}) \\
& = & {\rm gr}_{c}(H_{m})
\end{array}
$$
and since any finitely generated residually nilpotent group is Hopfian (see, \cite[Theorem 5.5]{mks}), we get the restriction of $\zeta_{m}$, say $\zeta_{m,c}$, on ${\rm gr}_{c}(H_{m})$ is a $\mathbb{Z}$-module isomorphism for all $c$ and $m \in \{2, \ldots, n\}$. For all $c$, let $\mu_{c}$ be the map from ${\rm gr}_c(H_2)\oplus\ldots\oplus {\rm gr}_c(H_n)$ into ${\mathcal L}^{c}_{1}(I_{n})$ defined by $\mu_{c}((u_{2}\gamma_{c+1}(H_{2}), \ldots, u_{n}\gamma_{c+1}(H_{n})) = \zeta_{2,c}(u_{2}{\rm I}_{c+2}{\rm A}(F_{n})) + \cdots + \zeta_{n,c}(u_{n}{\rm I}_{c+2}{\rm A}(F_{n}))$ for $u_{m} \in  \gamma_{c}(H_{m})$ for all $m \in \{2, \ldots, n\}$. It is easily checked that $\mu_{c}$ is a $\mathbb{Z}$-module isomorphism. By a result of Falk and Randell \cite[Theorem (3.1)]{fr}, ${\rm gr}_c(I_n)=\bigoplus_{i=2}^{n}{\rm gr}_c(H_i)$ and so, ${\rm gr}_{c}(I_{n})$ is isomorphic to ${\mathcal L}_1^{c}(I_{n})$ (as $\Z$-modules) by means of $\mu_{c}$ for all $c$. Therefore, ${\rm rank}({\mathcal L}^{c}_{1}(I_{n})) = {\rm rank}({\rm gr}_{c}(I_{n})) = {\rm rank}(L^{c})^{*}$. By Theorem \ref{3}, $L/J \cong {\rm gr}(I_{n})$ by a Lie algebra isomorphism $\overline{\xi}$ satisfying the conditions $\overline{\xi}(y_{ij}+J) = y_{ij}{\rm I}^{\prime}_n$ for all $2 \leq i \leq n$ and $1 \leq j \leq i$. Let $\zeta$ be the map from $L$ into ${\mathcal L}_{1}(I_{n})$ such that $\zeta(y_{ij}) = y_{ij}{\rm I}_{3}{\rm A}(F_{n})$ for all $2 \leq i \leq n$ and $1 \leq j \leq i$. Since $L$ is a free Lie algebra on the set $\{y_{ij}: 2 \leq i \leq n, 1 \leq j \leq i\}$, the map $\zeta$ extends uniquely to a Lie algebra epimorphism and so, $L/{\rm ker}\zeta \cong {\mathcal L}_{1}(I_{n})$. Clearly, $J \subseteq {\rm ker}\zeta$. Thus, $\zeta$ induces a Lie algebra epimorphism $\overline{\zeta}: L/J \rightarrow {\mathcal L}_{1}(I_{n})$ such that $\overline{\zeta}(y_{ij}+J) = y_{ij}{\rm I}_3{\rm A}(F_{n})$ for all $2 \leq i \leq n$ and $1 \leq j \leq i$. Moreover, $\overline{\zeta}$ induces a $\mathbb{Z}$-linear mapping $\overline{\zeta}_{c}: (L^{c}+J)/J \rightarrow {\mathcal L}^{c}_{1}(I_{n})$. By using similar arguments as in the proof of Theorem \ref{3}, we have $\overline{\zeta}$ is an isomorphism. Thus, $\overline{\zeta}(\overline{\xi})^{-1}$ is the required Lie algebra isomorphism from ${\rm gr}(I_{n})$ onto ${\mathcal L}_{1}(I_{n})$. \qed

\begin{remark}\label{thu1}\upshape{
By the proof of Theorem \ref{kostas}, we obtain ${\mathcal L}_1^{c}(I_{n}) = {\mathcal L}_1^{c}(H_{2})\oplus \ldots \oplus {\mathcal L}^c(H_{n})$. Since ${\mathcal L}_1^{c}(H_{i}) \cong \gamma_{c}(H_{i})/\gamma_{c+1}(H_{i})$ for all $i = 2, \ldots, n$, we get ${\rm rank}({\mathcal L}_1(H_{i})) = \frac{1}{c}{\displaystyle\sum_{d|c}\mu(d)i^{\frac{c}{d}}}$. Thus, for any $c$, $\frac{1}{c}{\displaystyle \sum_{d|c}\mu(d)\{2^{\frac{c}{d}}+ \ldots + n^{\frac{c}{d}}\}} \leq {\rm rank}({\mathcal L}^{c+1}({\rm IA}(F_{n}))$.
}
\end{remark}

\end{document}